\documentclass[a4paper,12pt]{amsart}

\usepackage {amsmath, amssymb, amscd, mathrsfs, amsthm, stmaryrd, bbm, diagbox, enumerate, slashed, graphicx, xspace, color, transparent,comment, subcaption}
\usepackage[all, cmtip]{xy}
\usepackage[text={6.3in,8.8in},centering,letterpaper,dvips]{geometry}
\usepackage{tikz, tikz-cd}
\usepackage{graphics}
\usepackage{amscd}
\usepackage{mathtools} 
\usepackage{amstext} 
\usepackage{array}   
\usepackage{algorithm}
\usepackage{algpseudocode}

\usepackage{url, hyperref}

\makeatletter
\newcommand{\addresseshere}{%
  \enddoc@text\let\enddoc@text\relax
}
\makeatother

\newcommand{\R}{\ensuremath{\mathbb{R}}}

\newcommand{\bea}{\begin{eqnarray}}
\newcommand{\eea}{\end{eqnarray}}
\newcommand{\be}{\begin{equation}}
\newcommand{\ee}{\end{equation}}

\def\H{\mathbb{H}}

\newtheorem{thm}{Theorem}[section]

\newtheorem{prop}[thm]{Proposition}

\newtheorem{conj}[thm]{Conjecture}

\theoremstyle{definition}
\newtheorem{rem}[thm]{Remark}

\usepackage{tikz}
\usetikzlibrary{decorations.pathreplacing,calligraphy}

\title[Searching for ribbons with machine learning]{Searching for ribbons with machine learning}

\author[Gukov]{Sergei Gukov}
\address{Dublin Institute for Advanced Studies, 10 Burlington Rd, Dublin, Ireland.\\
California Institute of Technology, Pasadena, CA 91125, USA}
\email{gukov@stp.dias.ie}

\author[Halverson]{James Halverson}
\address{Department of Physics, Northeastern University\\
360 Huntington Avenue, Boston MA 02115, USA. The NSF AI Institute for Artificial Intelligence and Fundamental Interactions.}
\email{j.halverson@northeastern.edu}

\author[Manolescu]{Ciprian Manolescu}
\address{Department of Mathematics, Stanford University\\
450 Jane Stanford Way, Building 380, Stanford, CA 94305-2125, USA.}
\email{cm5@stanford.edu}

\author[Ruehle]{Fabian Ruehle}
\address{Department of Physics and Department Mathematics, Northeastern University\\
360 Huntington Avenue, Boston MA 02115, USA. The NSF AI Institute for Artificial Intelligence and Fundamental Interactions.}
\email{f.ruehle@northeastern.edu}

\begin{document}

\setcounter{tocdepth}{2}

\begin{abstract}
We apply Bayesian optimization and reinforcement learning to a problem in topology: the question of when a knot bounds a ribbon disk. This question is relevant in an approach to disproving the four-dimensional smooth Poincar\'e conjecture; using our programs, we rule out many  potential counterexamples to the conjecture. We also show that the programs are successful in detecting many ribbon knots in the range of up to 70 crossings.
\end{abstract}
\maketitle

\section{Introduction}
Machine learning has recently been used in mathematics to suggest conjectures that later can be proved by humans, see for example~\cite{Carifio:2017bov} in the context of toric geometry,~\cite{Brodie:2019dfx}  for line bundle cohomology in projective spaces,~\cite{Bies:2020gvf} for Brill-Noether theory, and~\cite{DV,DJ} for knot theory and representation theory.  In this work, we illustrate how it can be used in a different way: to produce definitive results in mathematics. It is often perceived that neural networks can only establish approximate results, that may hold with high probability, but nevertheless are not  $100\%$ certain. While this may be the case for some problems that involve for example supervised learning, other learning algorithms can be devised to search for a path between two configurations using a particular set of steps. If successful, these produce a certificate which then can be verified by a human and is completely rigorous in the mathematical sense. 

The key point here is that the desired instance (e.g., a path between two configurations with a combinatorially large set of intermediate steps) may lie beyond current capabilities of the existent paper-and-pencil techniques. This is where machine learning comes to the rescue. It has the ability to quickly search through many potential solutions and, more importantly, to improve the search based on the successful ``games'' it plays. We use the word ``games'' since the same types of algorithms and architectures can be employed to play complex board games, such as Go or Chess \cite{Silver:2018aaa}, where the goals and winning strategies are similar to those in math problems. 

In our case, this general principle is realized in a class of problems in low-dimensional topology. The problems we have in mind involve finding a sequence of specific moves that relate knot diagrams with particular properties that we explain next. A {\it knot} (more generally, a {\it link}) is a circle (resp. a collection of circles) embedded in 3-dimensional space. Its projection to a generic 2-dimensional plane, called a {\it planar diagram}, is a collection of line segments (strands) that cross over or under each other. For a given knot or link $K$, the choice of diagram is not unique because one can continuously deform $K$ in the 3-dimensional space and also choose projections to different 2-planes. Nevertheless, one can show that different projections of the same  link $K$ are always related by a sequence of only three basic local operations on the planar diagram, the so-called Reidemeister moves. The converse is also true: if two diagrams are related by a sequence of such moves, then they represent the same link.

\begin{figure}
{
\fontsize{9pt}{11pt}\selectfont
   \def\svgwidth{5in}
   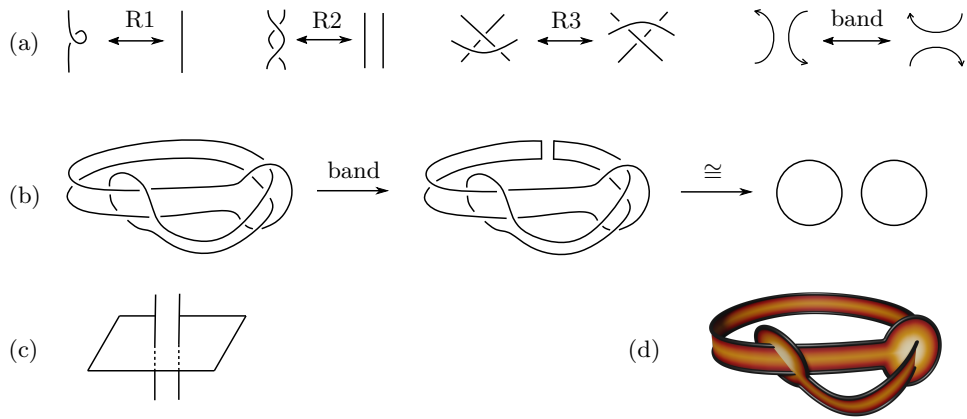
}
\caption{(a) The three Reidemeister moves and the band move. The band move needs to preserve orientations on the link. (b) After applying a band move to the square knot, the result can be deformed (via Reidemeister moves) into the unlink with $2$ components. (c) The kind of intersection allowed in a ribbon disk. (d) A ribbon disk for the square knot. }
\label{fig:Ribbon}
\end{figure}

Problems of this type (finding a sequence of moves between two diagrams) are ideally suited for machine learning because the set of all possible diagrams for a given link is huge (infinite, if we put no bound on the number of crossings) whereas the set of basic moves is very small. In other words, the search space is huge and there is no known algorithm that finds a solution other than brute-force; however, if someone gives you a set of moves that lead to a solution, they can be verified fast. Previously, reinforcement learning (RL) has been utilized for the unknotting problem: using Reidemeister moves to turn a complicated planar diagram into a trivial one \cite{Gukov:2020qaj}. In this paper, we tackle a much more difficult problem of the same type, where we extend the set of allowed operations on the knot diagram. Namely, to the set of the three Reidemeister moves we add the operation called {\em band addition}; these four moves are illustrated in Figure~\ref{fig:Ribbon} (a). Just like strands of a link can go either over or under each other, a band is allowed to go either over or under a strand and can pass through another piece of band.

The goal of our ``game'' then is to use these moves to transform a diagram of a given knot $K$ to that of a trivial link (a collection of $k$ split disjoint unknots), by using $k-1$ band moves (and any number of Reidemeister moves). Knots for which such sequences of moves exist are called {\it ribbon}. For example, the trefoil knot $3_1$ is not ribbon, but the square knot $3_1 \# (-3_1)$ pictured in Figure~\ref{fig:Ribbon} (b) is. There is currently no known algorithm for deciding whether a knot is ribbon. We use two different methods (Bayesian optimization and reinforcement learning) to search for moves that would prove that the knot is ribbon. While there is no guarantee that any ribbon knot will be proved to be ribbon in this way, our programs succeed in many cases. 

Detecting ribbons is relevant for a well-known strategy to disprove the smooth 4-dimensional Poincar\'e Conjecture (SPC4), a major open problem in topology. Counterexamples could be constructed if it can be shown that there exists a pair of knots with a certain common property (the same $0$-surgery), such that one knot is slice (bounds a disk in four dimensional half-space) and the other one is not. There are numerous topological knot invariants that are sufficient to establish that a knot is not slice, if the invariant does not take a particular value; these are known as slice obstructions. On the other hand, establishing that a knot is slice can be done by solving the ``ribbon game'' above: it is known that ribbon knots are slice. There are many examples of pairs of knots with common $0$-surgery; a large family of such examples was studied in \cite{MP21}. We used our techniques to show that for $843$ pairs in that family both knots are ribbon (and hence slice), thus eliminating  those potential counterexamples to SPC4.

In a different direction, we looked at knots with up to 14 crossings and established that 1705 of them are ribbon. (The same results were independently obtained by Dunfield and Gong \cite{DG}.)  We also developed two ways of generating ribbon knots of high crossing number, in the range of $15$ to $70$ crossings, and tested our algorithms on the resulting data.  These three data sets can be thought of as samples of ribbon knots from different data distributions, or priors on ribbon knots in the Bayesian context.

Surprisingly, in two different types of benchmarks on custom ribbon knot distributions, a Bayesian-optimized random walker outperformed RL. Still, RL outperformed a naive unoptimized random walker. There are a number of subtleties and caveats that are discussed in-depth in Section \ref{sec:Results} and summarized further in the conclusion.

The reader interested in finding ribbon disks for particular knots can experiment with our Bayesian-optimized random walker, which can be accessed at
\begin{center}
\href{https://github.com/ruehlef/ribbon}{https://github.com/ruehlef/ribbon}
\end{center}

\medskip
{\bf Acknowledgements.} We would like to thank Nathan Dunfield, Sherry Gong, Mark Hughes, and Lisa Piccirillo for helpful discussions during the preparation of this work. The code for attaching a band using the dual graph of the knot is based on previous work by Sherry Gong \cite{Gong}.

SG and CM are supported by a Simons Collaboration Grant on New Structures in Low-Dimensional Topology. CM is also supported by a Simons Investigator Award, and the NSF grant DMS-2003488. SG is also partially supported by the NSF grant DMS-1664227. JH and FR are supported by the National Science Foundation under Cooperative Agreement PHY-2019786 (The NSF AI Institute for Artificial Intelligence and Fundamental Interactions). JH is also supported by NSF CAREER grant PHY-1848089. FR is also supported by NSF grant PHY-2210333 and startup funding from Northeastern University.

\section{Ribbon disks} \label{sec:ribbon}
We can think of ribbon knots as those obtained from a trivial link by joining its $k$ components with $k-1$ bands. The result of the band addition operations is that a knot is ribbon if and only if it is the boundary of a {\em ribbon disk}, that is, one that lives in $3$-dimensional space and has self-intersections only of the form shown in Figure~\ref{fig:Ribbon} (c). 

If we allow ourselves one extra dimension---where the fourth coordinate is shown as the color in Figure~\ref{fig:Ribbon} (d)---then the self-intersections in the ribbon disk can be removed, and we obtain an embedded disk in $4$-dimensional space. Specifically, consider the half-space
$$ \H^4 = \{(x_0, x_1, x_2, x_3) \in \R^4 \mid x_0 \leq 0 \}$$
whose boundary is $\R^3$. A knot $K$ in $\R^3$ is called {\it slice} if there exists a smoothly embedded disk $D$ in $\H^4$ such that the boundary of $D$ is $K$.  We saw that every ribbon knot is slice. The converse is far from obvious and remains a famous unsolved problem in low-dimensional topology:

\begin{conj}[Slice-Ribbon Conjecture, \cite{MR0140100}]\label{conj:slice-ribbon}
Every slice knot is ribbon.
\end{conj}

Over the years, a variety of potential counterexamples to the Slice-Ribbon Conjecture have been constructed, usually involving knots with too many crossings to be analyzed by hand. For example, the authors of \cite{MR2740649} produce an infinite family of potential counterexamples, the simplest of which has 48 crossings.

The question of determining whether a given knot is slice (or ribbon) is of central importance in low dimensional topology. Indeed, it could help shed light on another major unsolved problem, namely the smooth Poincar\'e conjecture in dimension 4 (SPC4). This posits the non-existence of {\em exotic 4-dimensional spheres}, i.e., 4-dimensional smooth spaces (manifolds) that are topologically equivalent (homeomorphic) but not smoothly equivalent (diffeomorphic) to the standard 4-dimensional sphere  $S^4$.

\begin{conj}[SPC4]\label{conj:SPC4}
If a smooth 4-manifold is homeomorphic to $S^4$, then it is diffeomorphic to $S^4$.
\end{conj}

\begin{figure}[htb]
\centering
{
\fontsize{11pt}{11pt}\selectfont
   \def\svgwidth{3.3in}
   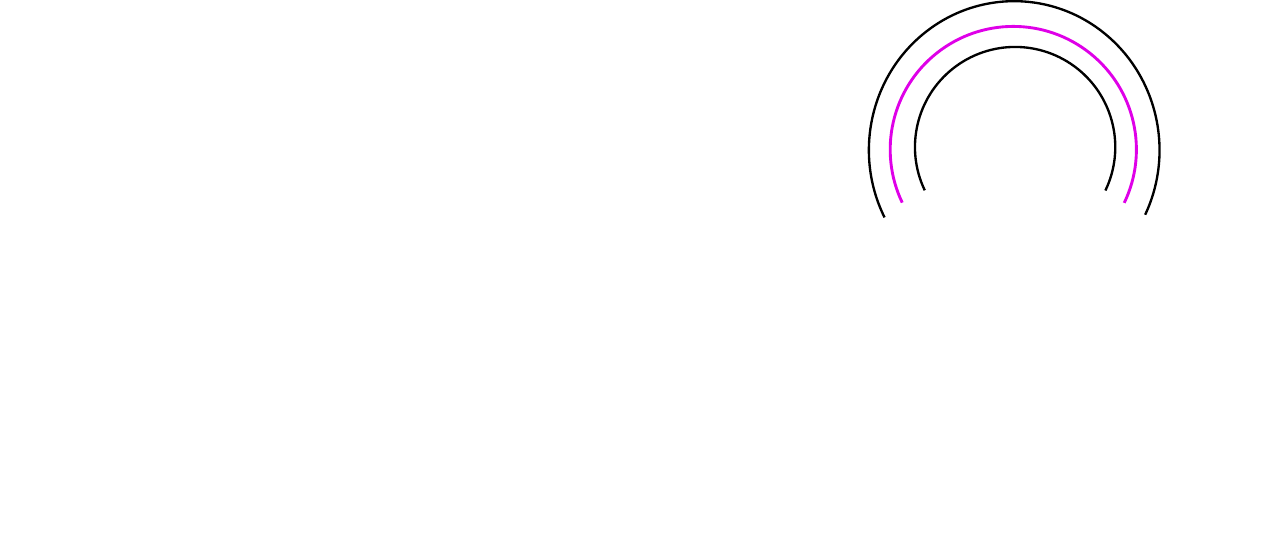
}
\caption{The exterior of a slice disk and the trace of the $0$-surgery. Left: We excise a tubular neighborhood of $\Delta$ (purple line) in $B^4$. Right: We attach a 2-handle to the boundary of $B^4$.}
\label{fig:Hsphere}
\end{figure}

Exotic spheres exist in many other dimensions, starting with dimension 7, and they are known not to exist in dimensions $1, 2, 3, 5$ and $6$. The four-dimensional problem remains open, and its relation to sliceness is due to the following source of potential counterexamples. Let $0$-surgery on a knot $K$ refer to the $3$-dimensional manifold $S^3_0(K)$ obtained from $S^3$ by removing a tubular neighborhood of the knot $K$, and gluing back a solid torus using a gluing map that swaps the meridian and the longitude:
$$ S^3_0(K) = (S^3 - \text{nbhd}(K)) \cup (S^1 \times D^2).$$
If one found a pair of knots which satisfy the following three properties:
\begin{enumerate}[(a)]
\item $K_1$ and $K_2$ have the same 0-surgery,
\item $K_1$ is slice,
\item $K_2$ is not slice,
\end{enumerate}
then an exotic $4$-dimensional sphere could be constructed. Indeed, let us view $B^4$ as obtained from $\H^4$ by attaching a point at infinity and consider the slice disk 
$$\Delta \subset \H^4 \subset B^4$$
with boundary $K_1 \subset S^3=\partial B^4$. After removing a standard neighborhood of $\Delta$ from $B^4$ we obtain a four-manifold $E(\Delta)$, called the {\em disk exterior},  whose boundary can be checked to be $S^3_0(K_1)$. On the other hand, one can produce another four-manifold $X(K_2)$, called the {\em trace} of the $0$-surgery on $K_2$, by starting with a ball $B^4$ and attaching a 2-handle $D^2 \times D^2$ (a neighborhood of a disk $D^2$, which is also topologically a 4-ball) to its boundary, where the attaching is done by gluing $\partial D^2 \times D^2$ to a tubular neighborhood of the knot $K_2$. The boundary of $X(K_2)$ is again $S^3_0(K_2)=S^3_0(K_1)$. See Figure~\ref{fig:Hsphere}.

 By gluing the two four-manifolds we just constructed along their common boundary, we obtain a closed $4$-manifold
$$ W = E(\Delta) \cup (-X(K_2)).$$
(The minus sign refers to a needed change in orientation.)
If we had considered $E(\Delta) \cup (-X(K_1))$ instead, from Figure~\ref{fig:Hsphere} we see that we would have gotten back $B^4 \cup B^4 = S^4$. In the situation at hand, one can still prove that $W$ is homeomorphic to $S^4$. Furthermore, by construction, the knot $K_2$ bounds a disk in $W \setminus B^4$. Since $K_2$ is not slice (i.e., does not bound a disk in $B^4 = S^4 \setminus B^4)$, we deduce that $W$ is not diffeomorphic to $S^4$ and is thus an exotic $4$-sphere.  We refer to \cite{MP21} for more details about this construction.

To pursue this strategy of disproving SPC4, one can search for potential candidates by first constructing pairs $(K_1, K_2)$ that satisfy condition (a) above. A systematic construction of all such pairs --- based on a certain class of 3-component links called ``RBG links'' --- was recently given in \cite{MP21}. Then, one needs obstructions to sliceness to further select pairs where one of the knots satisfies condition (c). Finally, one needs a method to show that the other knot in such a pair is slice, so that (b) is satisfied. In particular, since ribbon knots are slice, it suffices to show $K_1$ is ribbon. 

With regard to obstructions, topologists are able to show that certain knots are not slice (and hence not ribbon) using various invariants (numbers, polynomials, etc.) associated to the knots.  A simple example is the knot signature $\sigma (K)$ \cite{MR171275}, which must be $0$ for slice knots. Another is the Fox-Milnor condition on the Alexander polynomial $\Delta_K (x)$ of a knot $K$: if $K$ is slice, then this polynomial can be written in the form $\Delta_K (x) = f (x) f(x^{-1})$, where $f (x)$ has integer coefficients \cite{MR211392}. In practice, the conditions on the signature and the Alexander polynomial are already sufficient for most small knots; e.g. they can tell us that 2802 (94\%) of the 2977 prime knots with up to 12 crossings are not slice \cite{KnotInfo}. Moreover, the signature and the Alexander polynomial can be computed reasonably fast compared to other knot invariants. Nevertheless, it may happen that both conditions are satisfied for a non-slice knot, and yet one can prove it is not slice by using more sophisticated invariants and obstructions. Such obstructions come for example from knot homologies; see \cite{Rasmussen}, \cite{OStau}, \cite{Pic}, \cite{cablefig8}. Using all known obstructions, one can show that 17 more of the prime knots with up to 12 crossings are not slice, for a total of 2819 (95\%) out of 2977; see \cite{HKL}.

On the other hand, to show that a knot is ribbon, one typically wants to exhibit a ribbon disk; or, equivalently, a sequence of Reidemeister and band moves that take the knot into a trivial link. The remaining 158 prime knots with up to 12 crossings were shown to be ribbon using either paper and pencil or computer searches; see  \cite{Seeliger:2014aaa} and \cite{Lamm}. For larger knots, a computer program for this task was developed by Dunfield and Gong \cite{DG}. In the current work, building on their work, we use machine learning to expand the capabilities of the search for ribbon disks. 

For future reference, let us mention that there also exists an indirect way of proving that certain knots are ribbon, using the following result:
\begin{prop}[Gompf-Scharlemann-Thompson \cite{Gompf:2010aaa}]
\label{prop:GST}
If a two-component R-link consists of the unknot and another knot $K$, then $K$ is ribbon. 
\end{prop}

Here, a two-component link $L$ is called an {\em R-link} if $0$-surgery on it produces the connected sum $(S^1 \times S^2) \# (S^1 \times S^2)$. In \cite[Proposition 3.2]{Gompf:2010aaa} it is proved that if one component of such a link $L$ is the unknot, then $L$ can be transformed into the unlink by a sequence of handle slides. This immediately implies that the other component $K$ is ribbon; compare \cite[Section 8]{Gompf:2010aaa}.

\section{Generation of ribbon knots}
\label{sec:Generation}
As mentioned in the previous section, all prime knots up to 12 crossings have been completely classified as ribbon or non-ribbon. For prime knots up to 14 crossings, an almost complete classification of slice knots (for all but 21 knots) was obtained by Dunfield and Gong via brute force search of knots with no known slice obstruction; later, they expanded their classification efforts to knots with up to 19 crossings~\cite{DG}. Furthermore, Owens and Swenton~\cite{Owens:2021aaa} use a technique specific to alternating knots to identify ribbon disks in alternating knots up to 20 crossings. 

Since we want to study (ribbon) knots of a priori arbitrary crossing numbers, we implement two constructions which we will describe now. 

\begin{figure}
\centering
\includegraphics[width=5.5in]{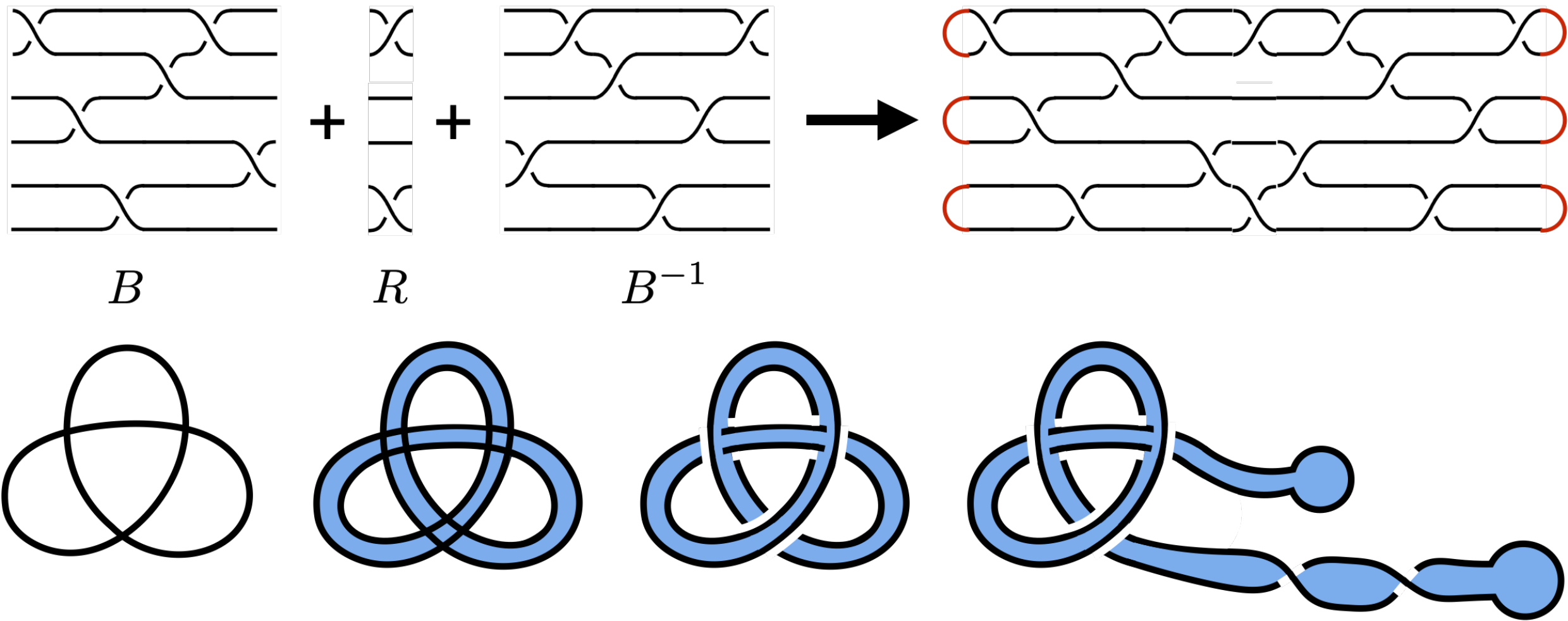}
\caption{Two ways of generating ribbon knots: ``Sym'' knots are obtained from symmetric reflections, which we generate using braids (top). ``Unsym'' knots are generated from doubling a random knot, choosing ribbon singularities, cutting it open, and inserting unknots on either end (bottom).}
\label{fig:KnotGeneration}
\end{figure}

The first construction, referred to as ``Sym'', uses the symmetric unions from~\cite{KT} which produce ribbon knots using the connected sum $K\#-\!K$ of a knot $K$ and its mirror, potentially with additional half-twists on the symmetry axis; see also \cite{Lamm}, \cite{Seeliger:2014aaa}. We generate such knots from random\footnote{Unless stated otherwise, random means drawn from a uniform distribution.} braid words $B$ (with generators $\sigma_i^{\pm1}$, an even number of strands $2m$ between 2 and 16, and a random number of generators chosen within a range such that the final knot has the desired number of crossings), inverting each generator to obtain $B^{-1}$, and building a random composition $R$ of $m-1$ non-consecutive generators $\sigma_i^{\pm1}$ to obtain the final braid word $W=B\circ R\circ B^{-1}$ for a symmetric knot $K$. The knot $K$ is generated from $W$ not by taking the usual braid closure (which would just result in a link with braid word $R$), but by taking the plat closure --- closing strands on either end of the braid word pairwise. The procedure is outlined at the top of Figure~\ref{fig:KnotGeneration}. The resulting link is simplified using SnapPy~\cite{SnapPy} and discarded if it has more than one component or if the final number of crossings is not in the desired range (simple rejection sampling).

\begin{figure}
\centering
\includegraphics[width=4.5in]{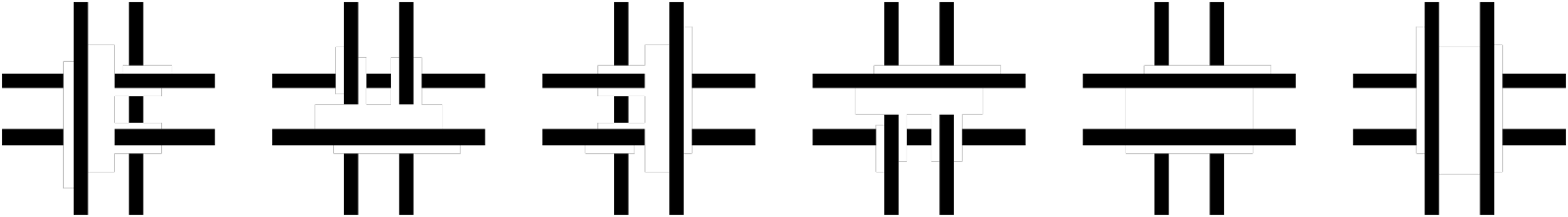}
\caption{Six possibilities for band intersections in a diagram.}
\label{fig:localChoices}
\end{figure}

The second generation method for ribbon knots, referred to as ``Unsym'', generates a random knot using SnapPy, with a random number of crossings again chosen within a range to produce a final knot with the desired crossing number. The procedure is outlined at the bottom of Figure~\ref{fig:KnotGeneration}: We start from the shadow of the random knot and double it to obtain a band. At each band intersection, we chose a random local picture from the six possible choices shown in Figure~\ref{fig:localChoices}. Subsequently, the band is cut open and we randomly insert up to two twists. Two unknots are inserted at the endpoints of the doubled tangle to obtain the ribbon knot. The resulting knot is simplified with SnapPy.

Apart from the data sets generated from the Sym and Unsym methods, we will also test our programs on the list of the $1705$ known prime ribbon knots with up to 14 crossings, following \cite{DG}. We will refer to this data set as ``Ribbon-to-14''.

\section{Machine learning}
To obtain verifiable truth certificates for a ribbon knot $K$, we identify a sequence of band insertions that produce an unlink from $K$ using reinforcement learning (RL) and Bayesian optimization of the environment associated to a Markov decision process.

\subsection{Markov Decision Process: The Band Environment}
\label{sec:BandEnvironment}
Our learning algorithms take place in the context of a {\em Markov decision process}, or environment, which requires defining a state space $\mathcal{S}$, an action space $\mathcal{A}$, and a state-dependent reward function $R$. 

Actions are chosen by sampling a policy function
\[
\pi: \mathcal{S} \to \mathcal{A}
\]
that is a probability density over actions. In general the policy $\pi$ is state-dependent, though we will also study cases in which it is state-independent up to a mild state dependence associated to masking out illegal actions. A sequence of actions is drawn from the policy, which determines a trajectory through state space (a game) that is intrinsically stochastic. Accordingly, the reward $R$ varies from one game to the next and is a stochastic function. A learning algorithm in this context aims to change the policy to optimize the expectation value of the reward. 

In the Markov decision process that we study, the action space $\mathcal{A}$ is given by
\begin{itemize}
\item \textbf{start:} Chooses an arc of the knot and starts a new band
\item \textbf{over:} Chooses the next arc and goes over this arc
\item \textbf{under:} Chooses the next arc and goes under this arc
\item \textbf{twist:} Inserts a twist into the band (either positive or negative)
\item \textbf{end:} Ends the band by attaching it to an arc. This completes a band addition move. We then simplify the resulting link using SnapPy, which includes all Reidemeister moves.
\end{itemize}
Illegal actions (which we mask out) consist of inserting positive and negative twists into the same band, self-intersecting the band, passing the band over or under itself, passing it over or under the same arc more than once, attaching bands that do not preserve orientation or connect different link components, and starting a new band before both ends of the current band have been attached. This allows us to model the band as a self-avoiding walk on the dual graph of the knot projection. We summarize the algorithm in Algorithm~\ref{alg:MDP} in Appendix~\ref{app:Algorithms}.

The states are specified by the current link plus the current position of the band. We describe this using four channels called $G,C,B,T$. These channels are $N\times N$ matrices where $N=c+2$ is the maximum number of nodes of the dual graph of a knot diagram with $c$ crossings. (Since band additions may increase the number of crossings, we will actually choose $N$ to be larger than $c+2$ of the initial knot; see the end of Section~\ref{sec:RL}.) The $G$ channel is a modification of the adjacency matrix of the dual graph of the current link, involving the crossing signs and vertex labels, so that we can uniquely reconstruct the link from this matrix. The $C$ channel is an (in theory redundant, but useful to the algorithm) channel that describes which nodes in the dual graph correspond to the same link component. The $B$ channel keeps track of how the band was routed by specifying a sequence of visited nodes in the dual graph. Finally, $T$ is just a constant matrix specifying the number of twists performed. 

\begin{figure}
    \centering
    \includegraphics[width=\textwidth]{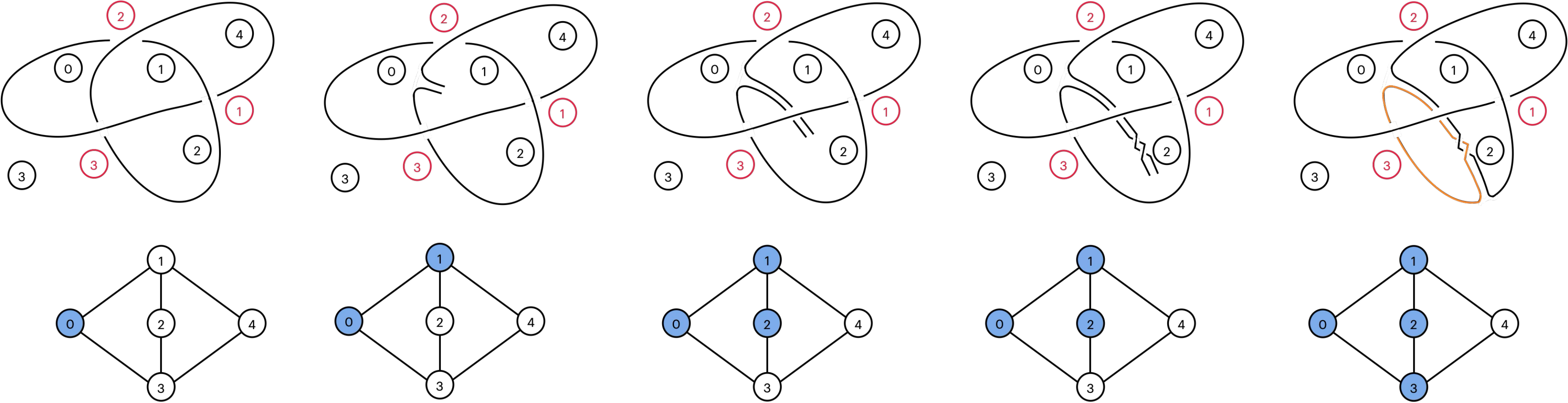}
    \caption{Adding a band to the trefoil and the corresponding dual graph.}
    \label{fig:GCBT-Example-Trefoil}
\end{figure}

We illustrate how these matrices are constructed and updated using a simple band in the trefoil knot, see Figure~\ref{fig:GCBT-Example-Trefoil}, the simplest non-trivial knot with $c=3$ crossings, and thus $N=5$. We give different stages of the knot as the band is formed in the top row, and the corresponding dual graph in the bottom row. Nodes in the dual graph (i.e., regions delineated by knot strands in the knot diagram) that have been visited are shaded blue in the dual graph. We have labeled the regions with numbers 0 to 4 as they appear in the dual graph. We have also labeled the crossings $c$ from one to 3 (in red). These numbers are used in the Graph matrix $G$.

To begin with we have no band (and hence no twists), as well as a one-component link. This means that the $B$ and $T$ matrices are simply $N\times N$ zero matrices,
\begin{align*}
    B =
\begin{bmatrix}
0 & 0 & 0 & 0 & 0 \\
0 & 0 & 0 & 0 & 0 \\
0 & 0 & 0 & 0 & 0 \\
0 & 0 & 0 & 0 & 0 \\
0 & 0 & 0 & 0 & 0
\end{bmatrix}\,,
\qquad
T =
\begin{bmatrix}
0 & 0 & 0 & 0 & 0 \\
0 & 0 & 0 & 0 & 0 \\
0 & 0 & 0 & 0 & 0 \\
0 & 0 & 0 & 0 & 0 \\
0 & 0 & 0 & 0 & 0
\end{bmatrix}\,.
\end{align*}
The component matrix $C$ is just the adjacency matrix of the dual graph, and the graph matrix $G$ is like the adjacency matrix, but modified to con tain signs and vertex labels that allow to uniquely reconstruct the link from the matrix,
\begin{align*}
C =
\begin{bmatrix}
0 & 1 & 0 & 1 & 0 \\
1 & 0 & 1 & 0 & 1 \\
0 & 1 & 0 & 1 & 0 \\
1 & 0 & 1 & 0 & 1 \\
0 & 1 & 0 & 1 & 0
\end{bmatrix}\,,
\qquad
G =
\begin{bmatrix}
0 & 3 & 0 & 2 & 0 \\
-3 & 0 & -3 & 0 & -1 \\
0 & 3 & 0 & 3 & 0 \\
-2 & 0 & -3 & 0 & -2 \\
0 & 1 & 0 & 2 & 0
\end{bmatrix}\,.
\end{align*}

The only legal first action is a start action, which means that we select a region adjacent to a strand on which we want to start the band. In this case we choose region 0 (see first diagram in Figure~\ref{fig:GCBT-Example-Trefoil}), indicated by adding a $1$ in the $(0,0)$-component of the band. The next legal action is an attach action, where we specify a second region adjacent to the strand. We choose region 1 in the example, so we set the $(1,1)$-component of $B$ to one. Together, the start and attach uniquely fix the strand on which we start the band, see second diagram of Figure~\ref{fig:GCBT-Example-Trefoil}). After the band is attached to a strand, we can move into a new region by crossing over or under one of the adjacent strands, or we can twist the band. In the example, we choose to cross under an adjacent strand into a new region, which we take to be region 3, see third diagram of Figure~\ref{fig:GCBT-Example-Trefoil}). We indicate that we moved from region 2 to region 3 by adding a $-1$ in the $(1,2)$ and $(2,1)$-component of the band matrix (the minus sign indicates an under-move). These actions change the band matrix according to
\begin{align*}
    B =
\begin{bmatrix}
0 & 0 & 0 & 0 & 0 \\
0 & 0 & 0 & 0 & 0 \\
0 & 0 & 0 & 0 & 0 \\
0 & 0 & 0 & 0 & 0 \\
0 & 0 & 0 & 0 & 0
\end{bmatrix}
\xrightarrow{\text{start}}
\begin{bmatrix}
1 & 0 & 0 & 0 & 0 \\
0 & 1 & 0 & 0 & 0 \\
0 & 0 & 0 & 0 & 0 \\
0 & 0 & 0 & 0 & 0 \\
0 & 0 & 0 & 0 & 0
\end{bmatrix}
\xrightarrow{\text{attach}}
\begin{bmatrix}
1 & 0 & 0 & 0 & 0 \\
0 & 1 & 0 & 0 & 0 \\
0 & 0 & 0 & 0 & 0 \\
0 & 0 & 0 & 0 & 0 \\
0 & 0 & 0 & 0 & 0
\end{bmatrix}
\xrightarrow{\text{under}}
\begin{bmatrix}
1 & 0 & 0 & 0 & 0 \\
0 & 1 & -1 & 0 & 0 \\
0 & -1 & 1 & 0 & 0 \\
0 & 0 & 0 & 0 & 0 \\
0 & 0 & 0 & 0 & 0
\end{bmatrix}\,.
\end{align*}
After these steps, we decide to twist the band clockwise twice, see diagram 4 of Figure~\ref{fig:GCBT-Example-Trefoil}). This does not change the $B$-matrix, but is reflected in adding $-1$ to the $T$ matrix for each twist (minus because we twisted clockwise rather than counter-clockwise),
\begin{align*}
    T =
\begin{bmatrix}
0 & 0 & 0 & 0 & 0 \\
0 & 0 & 0 & 0 & 0 \\
0 & 0 & 0 & 0 & 0 \\
0 & 0 & 0 & 0 & 0 \\
0 & 0 & 0 & 0 & 0
\end{bmatrix}
\xrightarrow{\text{twist}}
\begin{bmatrix}
-1 & -1 & -1 & -1 & -1 \\
-1 & -1 & -1 & -1 & -1 \\
-1 & -1 & -1 & -1 & -1 \\
-1 & -1 & -1 & -1 & -1 \\
-1 & -1 & -1 & -1 & -1
\end{bmatrix}
\xrightarrow{\text{twist}}
\begin{bmatrix}
-2 & -2 & -2 & -2 & -2 \\
-2 & -2 & -2 & -2 & -2 \\
-2 & -2 & -2 & -2 & -2 \\
-2 & -2 & -2 & -2 & -2 \\
-2 & -2 & -2 & -2 & -2
\end{bmatrix}\,.
\end{align*}
As the next action, we choose an attach action which ends the band. We choose to end it in region 3, see diagram 5 of Figure~\ref{fig:GCBT-Example-Trefoil}). We do not update the band matrix to reflect this change; rather, since the band is now completed, we reset it (and the twist matrix) to the zero matrix, in preparation for the next band. But we do update the graph and component matrix to describe the new link that we obtained from attaching this band. We highlighted the two link components in black and orange. The resulting two-component link consists of a trefoil (black) and a linked unknot (orange). This is a legal band, and after simplification, which includes Reidemeister move 2 on the orange component, we would have added one new crossing, the twist in the band. 
\begin{rem}
\label{rem:more}
Our representation of the link in terms of the GCBT channels does not capture self-intersecting bands or the possibility that the same node on the dual graph is visited more than once. Mathematically, such bands are allowed (and perhaps even necessary in a given diagram). These more general bands can be simulated in our set-up by making the diagram larger using Reidemeister moves. Including them in the action space associated to a fixed diagram would require a fundamentally different data structure, which might be worthwhile exploring in the future.
\end{rem}

In our set-up, terminal states that ``win the game'' are unlinks. Hyperparameters include a maximum number of steps, number of dual graph vertices, number of components, and number of bands. If during the game these hyperparameters are exceeded, or if the band cannot be ended consistently with the legal actions described above, the game is ``lost." An episode experienced by the agent goes from the start of the game until it is won or lost.

Since Bayesian optimization and reinforcement learning use different reward functions and methods of optimizing the Markov Decision Process, we introduce them in two different sections.

\subsection{Bayesian optimization}
We consider first the case that the policy $\pi$ is a state-independent probability distribution on actions that allows for the four different action types outlined above to be weighted differently. 

We wish to find optimal weights $w$ for gameplay. In this case, the reward function $R$ is the number of successfully recognized Sym or Unsym ribbon knots in a range from 15-50 crossings (with 100 knots per crossing) and where the game is played for up to 5~min per knot and we allow for up to 5 bands. If we knew the expectation value of this function as a function of the hyperparameters, we could optimize it directly using for example gradient-based optimizers. However, we do not know it, and each run is extremely time-expensive, since it involves scanning over $100\times 36=3600$ knots for potentially up to 5 minutes (worst case $300$ core hours). Hence, we opted to choose a Bayesian hyperparameter optimization scheme which requires few runs and models the stochastic reward function  based on previous results.

In Section~\ref{sec:BandEnvironment}, we listed 5 types of actions. In general, these depend on the current state (i.e., the link and the current position of the band) and in any given state, most actions are actually illegal. Nevertheless, we model the 5 action types to be state-independent aside from masking out all illegal actions in a given state. For example, once a band is started, starting another band is illegal until the current band has been attached. In contrast, for the over/under and attach moves, there are typically more than one possibility. In such a case, performing one of these actions with probability $p$ means that any of the possible legal moves for this action will be performed with possibility $p$.

To make the search more feasible, we further minimize the number of independent hyperparameters for the 5 actions: since the distributions of knots is symmetric with respect to going over or under an arc for any given band, we assign these two actions equal weights in the Bayesian optimization procedure. Also, since an overall scaling of all weights results in the same probability distribution, we fix the weight of one of the actions (over and under, say) to 1. Furthermore, since start actions are either all illegal (if a band is already started) or the only legal actions (if no new band has been started), it never competes with other action types and we may hence fix it to 1 as well. In addition to the actions, the maximum number of actions that the agent is allowed to perform before the episode counts as a loss is a very important hyperparameter which we also want to optimize. Hence, we have 3 hyperparameters: the two independent weights $w_{1,2}\in \mathbb{R}^2$ and the maximum number of actions per episode, $w_3\in\mathbb{R}$, which we take to be in $[10,10000]$. 
Then, for any value of these hyperparameters, one may play many games in the band environment and collect rewards, which together make up a set of observations $\{w_n, R(w_n)\}$. 

To optimize the expected value of the reward, we wish to model the reward via a stochastic process. A convenient choice that facilitates optimization via Bayesian inference is to model $R(w)$ as a draw from a Gaussian process (GP), meaning that for any finite set of weights $\{w_i \in \mathbb{R}^3\}$, the vector of random variables $R(w_i)$ is distributed according to a multivariate Gaussian. Following \cite{Snoek:2012aaa} and using the implementation of~\cite{BayesianOptimization:2023aaa}, we utilize Bayesian optimization to find weights that lead to better performance of the state-independent agents in the band environment. It is assumed that $R(w)$ is drawn from a GP prior with a chosen mean and covariance. Beginning with the prior, the algorithm iteratively chooses the next weight to observe $w_\text{next}$, plays a game at $w_\text{next}$ and records the reward, computes the Bayesian posterior given the new data point $\{w_\text{next}, R(w_\text{next})\}$, and records the best weight seen so far $w_\text{best}$. At each step of the optimization, the Bayesian posterior is Gaussian and depends on a collection of game results $\{w_n,R(w_n)\}$ with predictive mean function $\mu(w; \{w_n,R(w_n)\})$ and variance $\sigma^2(w;\{w_n,R(w_n)\})$. At any given step a useful quantity is
\[
\gamma(w) = \frac{R(w_\text{best})-\mu(w;\{w_n,R(w_n)\})}{\sigma(w;\{w_n,R(w_n)\})},
\]
which measures the deviation of the mean prediction at $w$ from the current best and normalizes by the variance. 
More specifically, the next weight is chosen as
\[
w_\text{next} = \text{argmax}_w \, a(w)
\]
according to an acquisition function $a:\mathbb{R}^3 \to \mathbb{R}^+$ that maximizes the expected improvement,
\[
a(w; \{w_n, R(w_n)\}) = \sigma(w;\{w_n,R(w_n)\})\left[\gamma(w)\Phi(\gamma(w)) + \mathcal{N}(\gamma(w); 0, 1)\right],
\]
where $\Phi(\cdot)$ is the cumulative distribution function of the standard normal. \cite{Snoek:2012aaa} found this acquisition function to be superior to two others they considered.

We run the optimizer with 10 initialization points and perform 50 iterations. For the optimal hyperparamters, the optimizer finds that it is beneficial to give a lot of weight to attaching as compared to continuing the band using over/under moves or twisting. This means that most knots, at least in the datasets we tested, can be shown to be ribbon using a collection of short bands. Note that this might mean that we are using more but simpler bands than the minimum fusion number would require. 

For the number of maximum steps taken before we force a reset to the original knot, we find that a large value gives the best results. Resetting earlier could be beneficial if the agent attached a sequence of bands that produce a link which cannot be ribbon. Adding more bands will not change this, so it might be better to start over. (Note that we do not check slice obstructions after band addition, since it is quite costly for larger links.) From looking at the steps taken, however, the situation described above seems to be rather rare: for ribbon knots there often seem to be several sequences of bands that will lead to the unlink. Moreover, the agent tends to get stuck, meaning there is not a single legal action left for it to do, after a few hundred steps, which automatically forces a reset. These two effects combined mean that as soon as the maximum steps hyperparameter exceeds a few hundred, the link is often either solved (meaning a set of ribbon bands is found) or reset to the starting link (meaning the agent got stuck). In either case, the precise value of the hyperparameter is not meaningful, and the Bayesian optimizer found a wide variety of values that seemed to work equally well. For concreteness, we chose to set it to 5500. 

For the relative weights of the actions, the optimum that is found by the Bayesian optimizer is
\begin{align}
   \label{eqn:weights bayes}
[\text{start\;:\;end\;:\;over\;:\;under\;:\;twist}] = [1: 17: 1 : 1 : 3]\,.
\end{align}

This means that each possible move of a given move type is given the associated weight, and the agent selects from the associated distribution on actions after masking out illegal moves.
We refer to a random walker with this ratio of different move types as a Bayes RW. On the other hand, we refer to the agent with ratios
\begin{align}
   \label{eqn:weights naive}
   [\text{start\;:\;end\;:\;over\;:\;under\;:\;twist}] = [1 : 1 : 1 : 1 : 1]\,.
\end{align}
as the Naive RW, since it gives equal weight each of the different action types.  For either RW, probabilities are sampled from the stated distributions with a canonical state-independent sampler, implemented in Python as \texttt{np.random.choice}.

In the context of the algorithm in Appendix~\ref{app:Algorithms}, this means that the policy $\pi(a_t~|~s_t)$ on line 6 is implemented as
\begin{align*}
    \pi(a_t~|~s_t) = \texttt{np.random.choice}([\text{start,end,over,under,twist}], p=[1,17,1, 1,3])\,,
\end{align*}where over represents $N$ actions going over any of the $N$ strands (and similarly for under), but we mask out all strands that cannot be crossed from the given strand position

\subsection{Reinforcement Learning} \label{sec:RL}
Reinforcement learning (RL) is another learning mechanism for optimizing a Markov decision process. It utilizes a state-dependent policy function that is approximated by a deep neural network, in which case the learning is known as deep RL. There are a variety of different RL algorithms, roughly classified into so-called policy-based methods, which directly optimize the policy, and value-based methods, which implicitly optimize the policy via optimizing the so-called value function of a state. The value function is
\begin{align*}
V(s) = \mathbb{E}[G_t| s_t = s],
\end{align*}
where $s$ is a state and
\begin{align*}
G_t = \sum_{k=0} \gamma^k \, R_{t+k+1}
\end{align*}
is an accumulated reward known as the return, which is discounted by a discount factor $\gamma \in (0,1]$. The value function measures the expected return across many trajectories through state space determined by draws from the policy. The action-value function 
\begin{align*}
   Q(s,a) = \mathbb{E}[G_t| s_t = s,~a_t=a],
\end{align*}
measures the value of pairs of states and actions. These functions are in general not known and are also estimated by neural networks. The advantage function 
$A(s,a) := Q(s,a) - V(s)$  measures the difference in value between conditioning on the initial state and action, as opposed to just the state. The RL algorithm then proceeds by playing episodes (i.e., running through the loop in Algorithm~\ref{alg:MDP}) multiple times, computing the return $G$ (or advantage), and then trains neural networks that approximate the policy $\pi(a|s)$, value function $V(s)$, and/or action-value function $Q(s,a)$. The policy NN is trained to maximize $G$, while the state value and action value function are trained to accurately predict $V$ and $Q$. See 
\cite{Halverson:2019tkf,Ruehle:2020jrk} for a more thorough introduction in the string theory literature,
and \cite{Gukov:2020qaj} for an application of RL to the unknotting problem. 

For our application of RL, we utilize the band environment described in Section~\ref{sec:BandEnvironment}. For the rewards we tried combinations of various different possibilities. In every case, we reward the agent by a fixed amount if it wins the game. In some cases, it receives an intermediate reward given by the change in the number of crossings of the link after attaching a band (note that the winning terminal state, the unknot, has zero crossings, such that the maximal reward for a knot with $c$ crossings is $c$), or it was punished for each step taken (in addition to the discount factor $\gamma<1$). In these cases, we trained the RL algorithm on ribbon knots of ascending difficulty, where we took the crossing number as a surrogate for the difficulty measure of finding bands and hence establishing ribonness.

For RL algorithms we tried TRPO~\cite{Schulman:2015TRPO} and A3C~\cite{Mnih:2016A3C}, which have  been applied successfully by the authors in the past to problems in knot theory~\cite{Gukov:2020qaj} and a number theory problem arising in string theory~\cite{Halverson:2019tkf}. We also implemented an adaptation of the AlphaZero Monte-Carlo tree search (MCTS) algorithm~\cite{Silver:2018aaa} to work with single-player games. 
TRPO is optimizing the policy function using a trust region and a line search. A3C uses asynchronous actor-critics to update the value and policy functions based on the advantage, introduced above. AlphaZero performs a Monte-Carlo tree search (MCTS) which consists of simulation steps where games are rolled out to a terminal state. 

From our experiments, we found that TRPO seemed to perform on the same level as A3C for this problem. This was somewhat surprising, since TRPO vastly outperformed A3C for the (similar but much simpler) unknot problem \cite{Gukov:2020qaj}. For the MCTS using the AlphaZero implementation, we were not able to obtain results since the simplification of the link in the Monte Carlo Tree Search is very costly and we lack the computing resources to train this algorithm for a meaningful number of steps.

Let us now describe in detail the RL methods utilized to produce the results in Section \ref{sec:Results}. The algorithm that we choose, for some of the reasons mentioned, is A3C. The value function is approximated by a neural network composed of a convolutional layer with 16 output channels, kernel size $6$, and stride $3$; a sigmoid non-linearity; a max-pool with stride and kernel size $2$; a linear layer of width $128$; a sigmoid non-linearity; and a linear layer of width $1$, the output of the neural network that models the scalar value function. This architecture is similar to LeNet. The neural network that approximates the policy is similar, but the output dimension of the last linear layer is $|\mathcal{A}|$, the size of the action space, and a Softmax is applied at the end to turn it into a probability density. The input to both neural networks consists of the GCBT channels for the last step. We also tried a ResNet architecture, but found that the sampling time of this larger network was too costly in comparison to the policy improvements in the timeout benchmarks described below.

A key point of our analysis is that the last layer of the policy network has weights and biases that are initialized to zero. Due to the Softmax, this means that at initialization the policy $\pi$ is a uniform distribution over the space of possible actions. We refer to this as a {\em Naive start} for the RL training, as it is equivalent to  Naive RW, aside from the fact that the distribution is sampled by network calls instead of \texttt{np.random.choice}; we will see the importance of this distinction. Alternatively, the uniform distribution associated to the policy network can be re-weighted to match \eqref{eqn:weights bayes}, in which case the initial policy network is equivalent to the Bayes RW. We refer to this as a {\em Bayes start} for RL training.

We scanned over many hyperparameters attempting to optimize the performance of the trained agents according to detailed benchmarks described below. In the end, the results we present arise from RL experiments with $\gamma = .9$, up to $10$ band additions,  and up to $10$ link components at intermediate stages. For Sym, Unsym, and Ribbon-to-14 knots we used a maximum number $N$ of dual graph nodes of $100$, $100$, and $20$, respectively. Evaluation during training and in post-training benchmarks is performed on a test set of $200$ knots for Sym and Unsym, out of a total database of $10000$ for each. For Ribbon-to-14, all $\sim1700$ knots are used for evaluation and benchmarks; since there is less data, we wish to train on all of it, not holding out a test set.

\subsection{Benchmarks}

Having discussed Bayesian Optimization and RL in our band environment, we now have six agents on which to run benchmarks: Naive RW and Bayes RW, which are defined by state-independent policies that are sampled with \texttt{np.random.choice}; untrained policy neural networks with Bayes or Naive starts, which are equivalent distributions to Bayes RW and Naive RW but sample actions instead from the neural network; and policy neural networks that are trained with A3C from a Bayes or Naive start. We perform two types of benchmarks for the agents:
\begin{itemize}
\item {\em number-of-episodes benchmarks} in which an agent is run (on a fixed set of knots) for up to a fixed number of episodes; 
\item {\em timeout benchmarks}, in which the agent is run for up to a fixed number of seconds.
\end{itemize}
The knots utilized in these post-training benchmarks are the same as those utilized in evaluation during training. An episode ends if there are no legal moves, or if there are too many dual graph vertices, bands, or components.

\section{Results}
\label{sec:Results}

\begin{figure}[t]
   \centering
   \includegraphics[width=\textwidth]{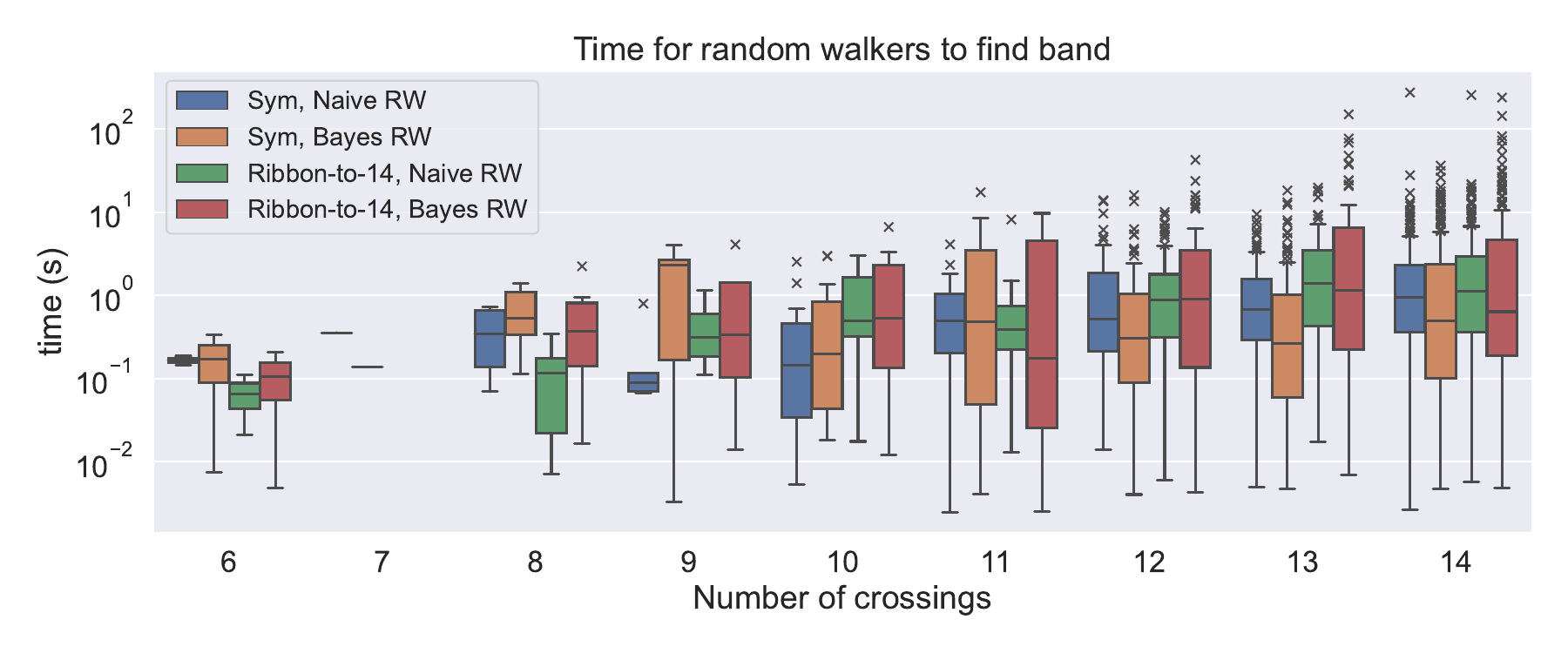}
   \caption{Performance of the Bayes RW and Naive RW on the Ribbon-to-14 dataset and ribbon knots in the same crossing range generated using the Sym generator.}
   \label{fig:sym_vs_rto14}
\end{figure}

\begin{figure}[t]
   \centering
   \includegraphics[width=.46\textwidth]{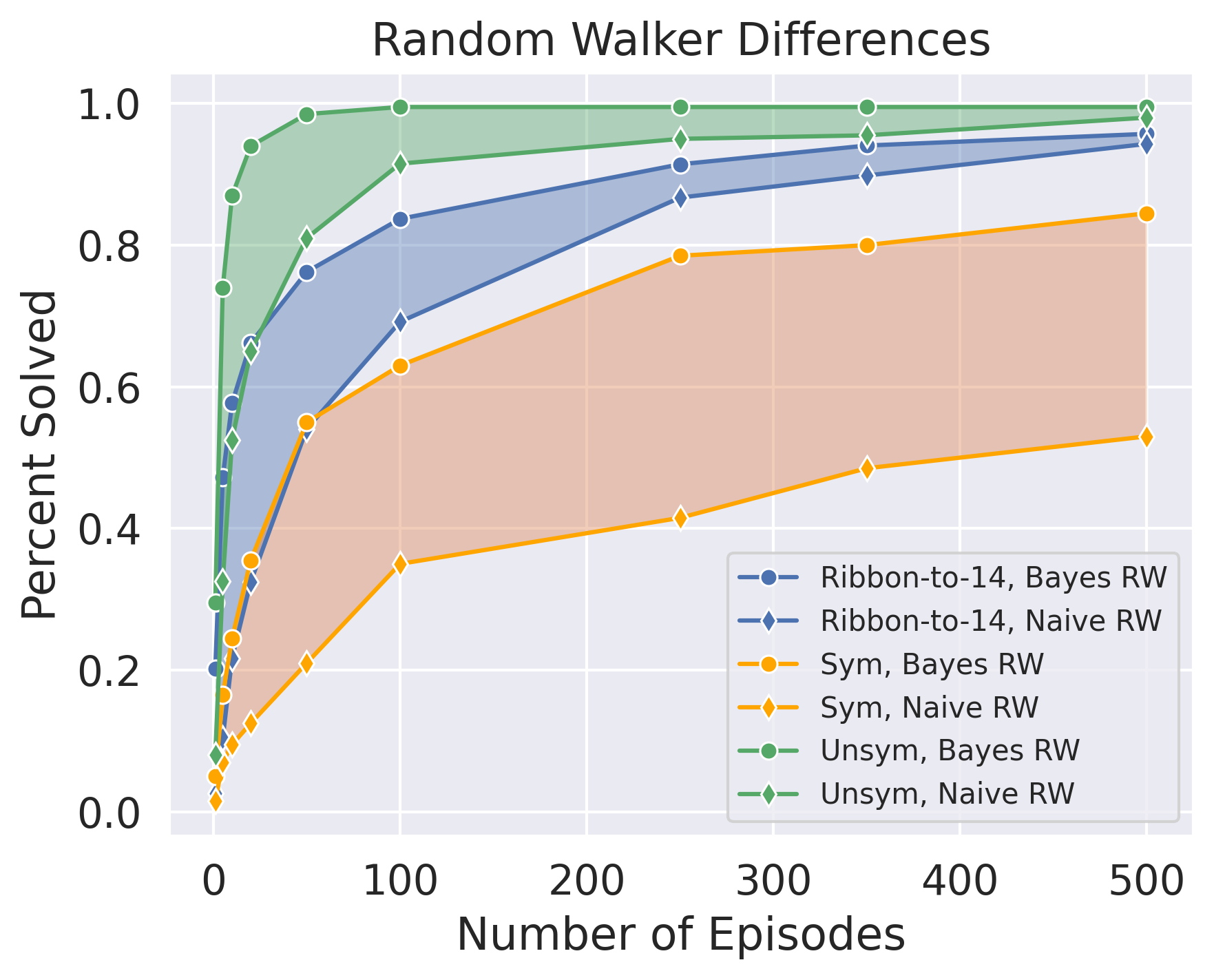}~
   \includegraphics[width=.46\textwidth]{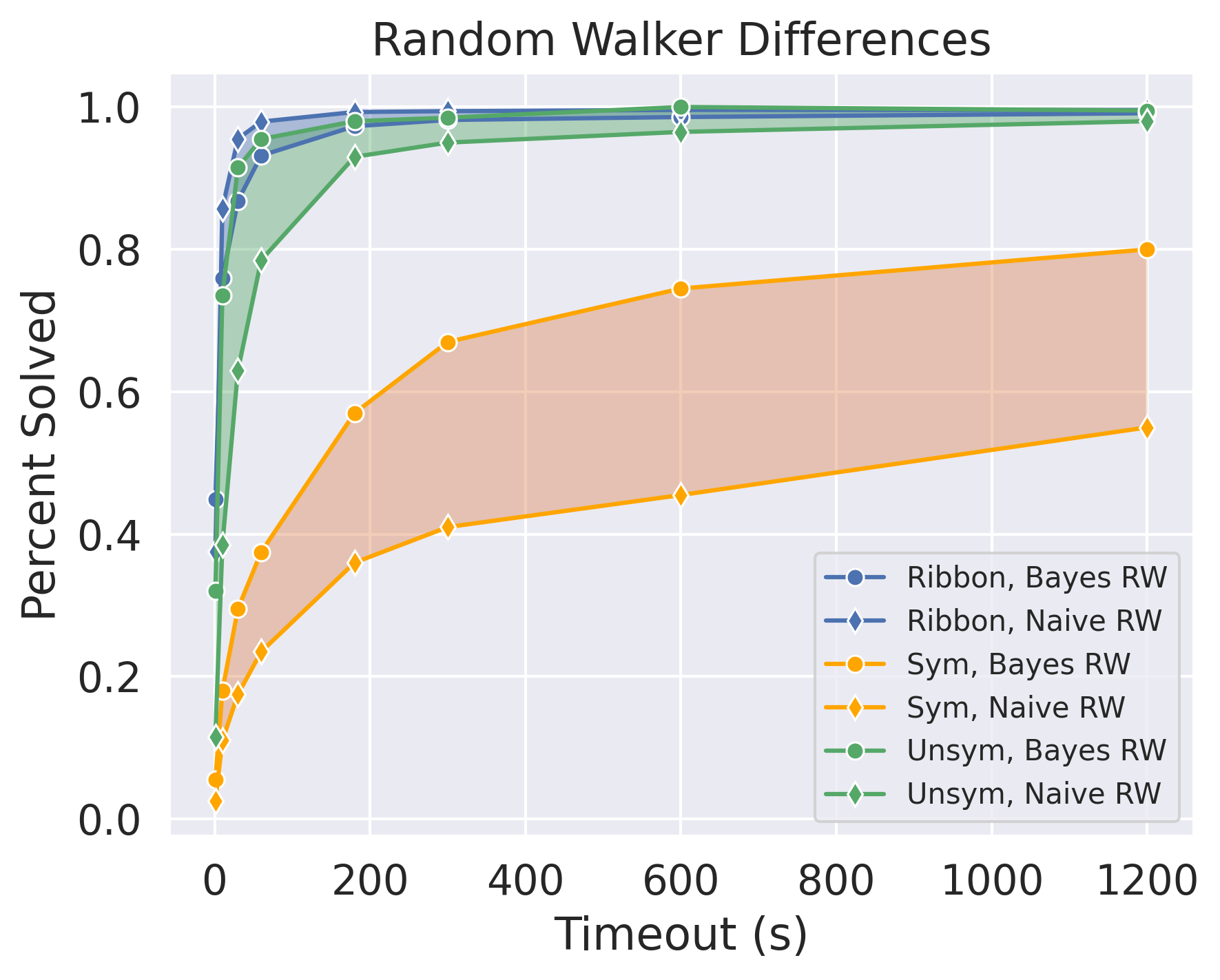}
   \caption{Performance differences between the Bayes RW and Naive RW for different types of knots. Naive RW performs better only for Ribbon-to-14 timeout benchmarks.}
   \label{fig:envelopes}
\end{figure}

\begin{figure}[t]
   \centering
   \includegraphics[width=.46\textwidth]{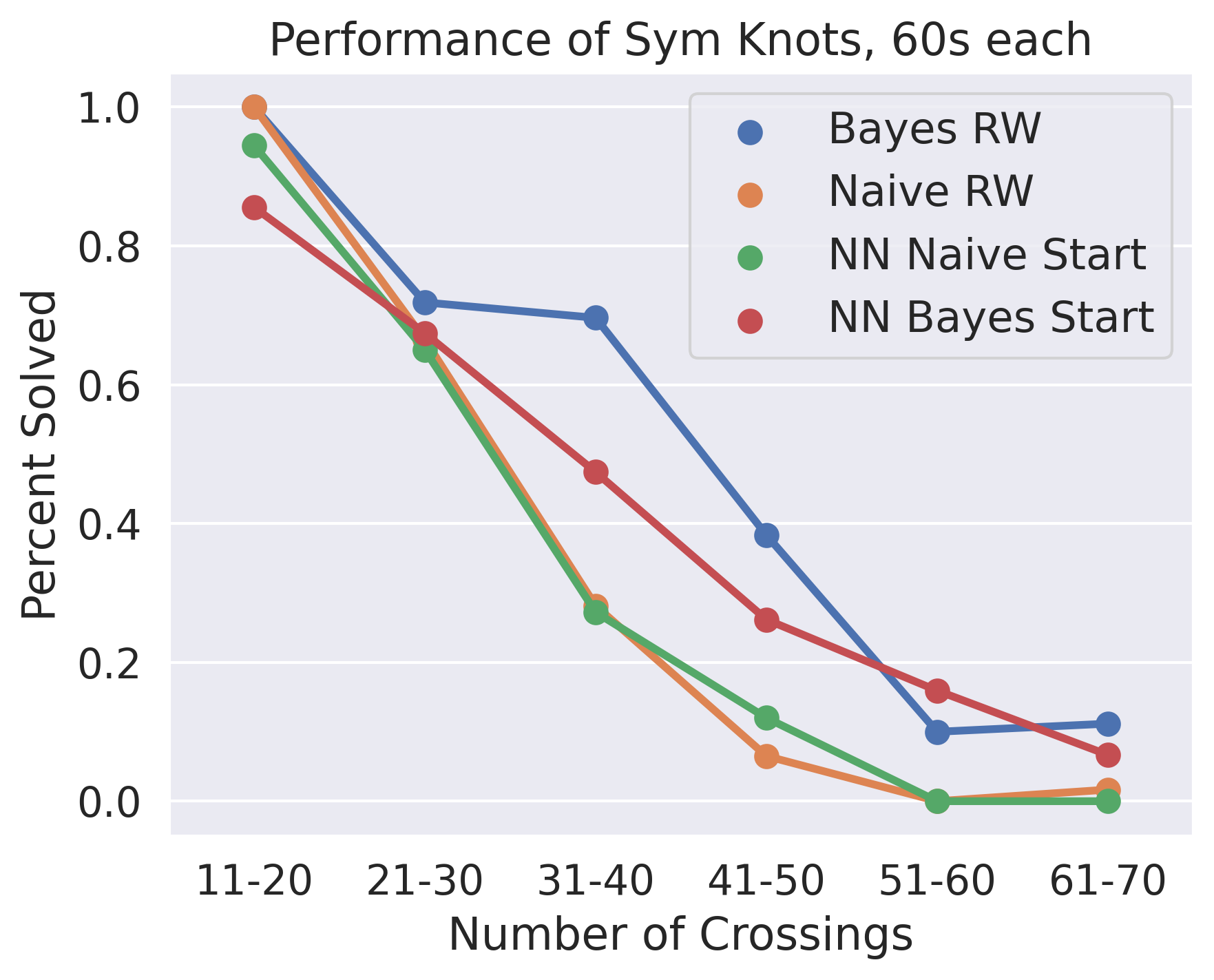}~
   \includegraphics[width=.46\textwidth]{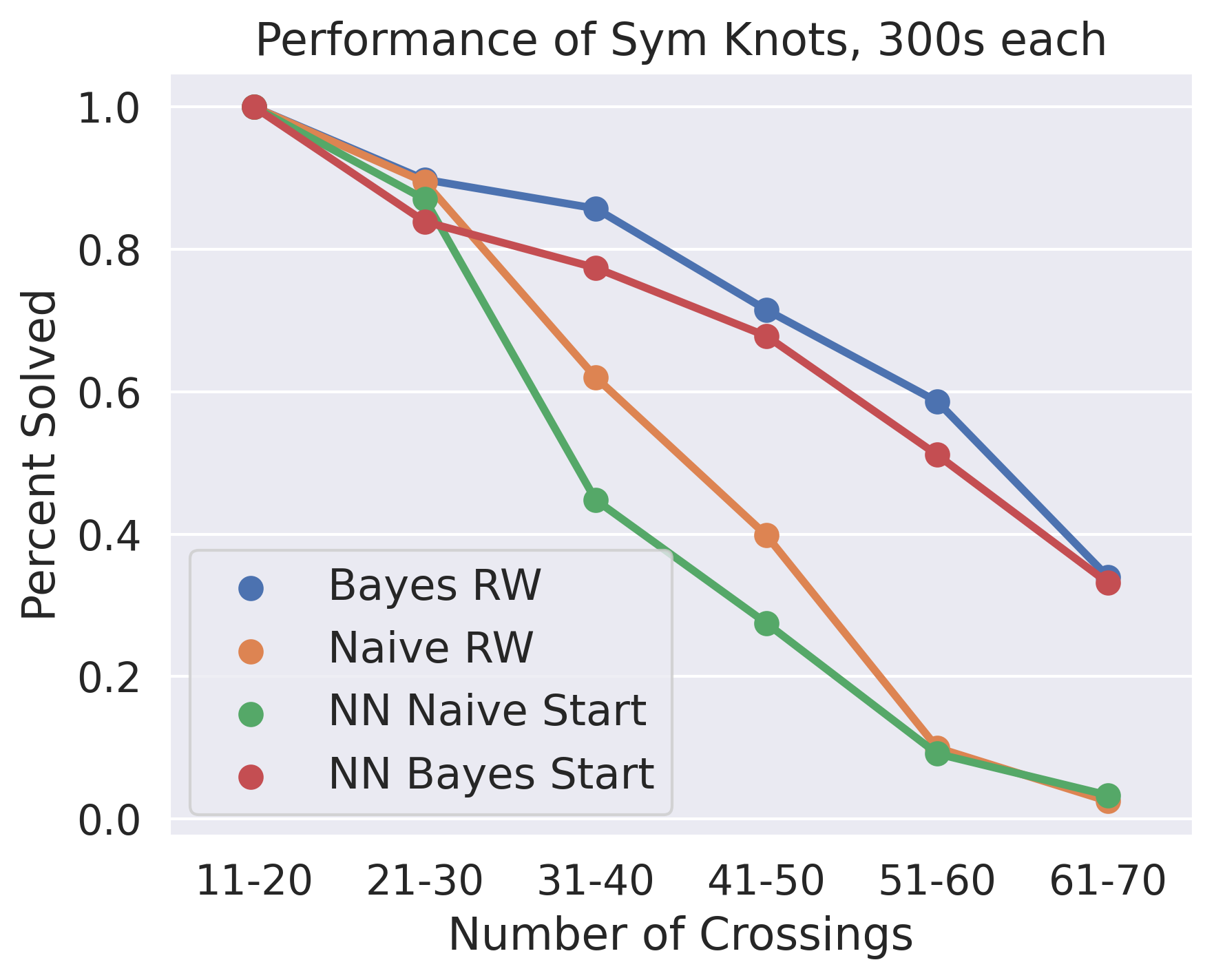}\\   
   \includegraphics[width=.46\textwidth]{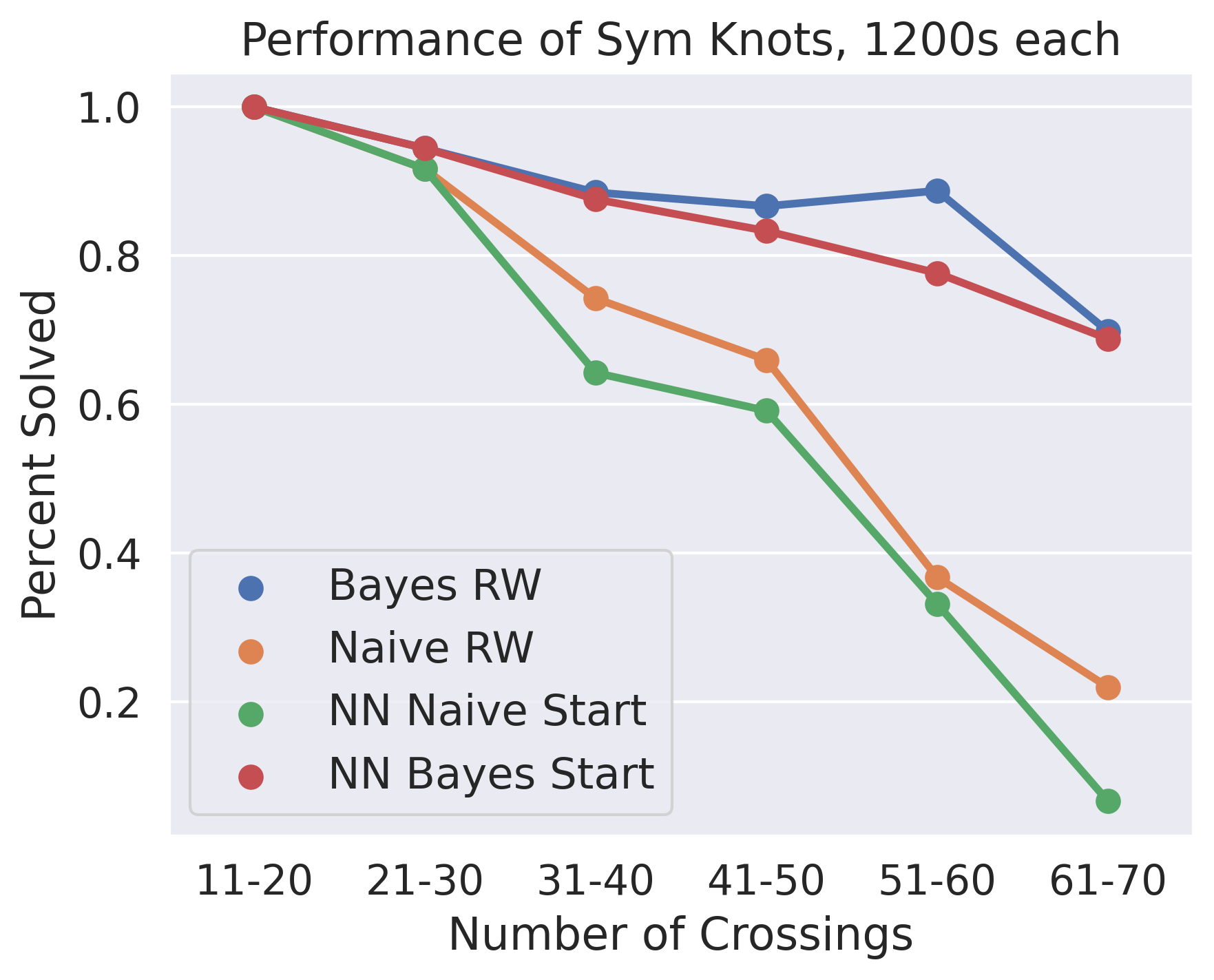}~
   \includegraphics[width=.46\textwidth]{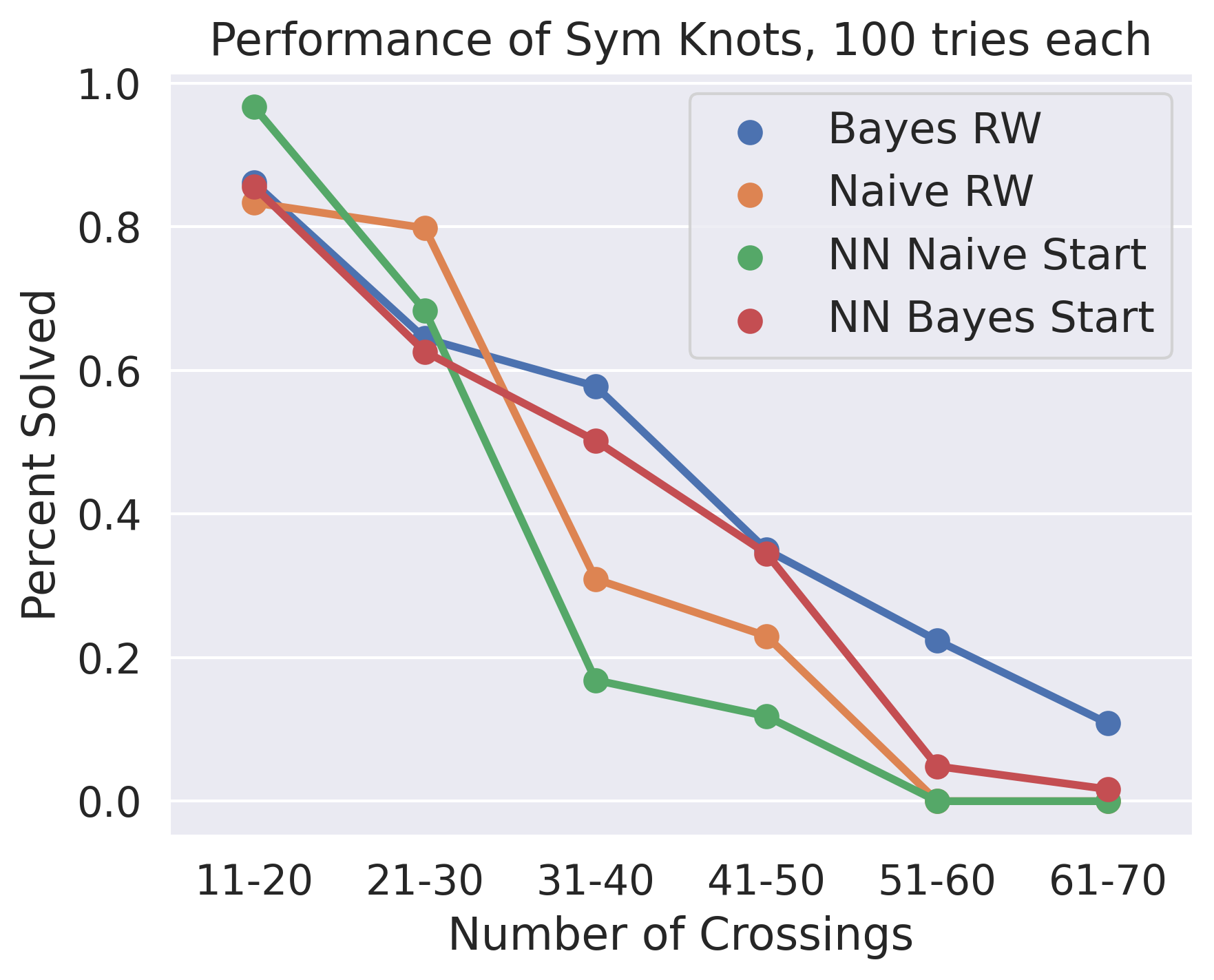}
   \caption{Dependence of performance on the number of crossings for Sym knots.}
   \label{fig:n_dep_sym}
\end{figure}

\begin{figure}[t]
   \centering
   \includegraphics[width=.46\textwidth]{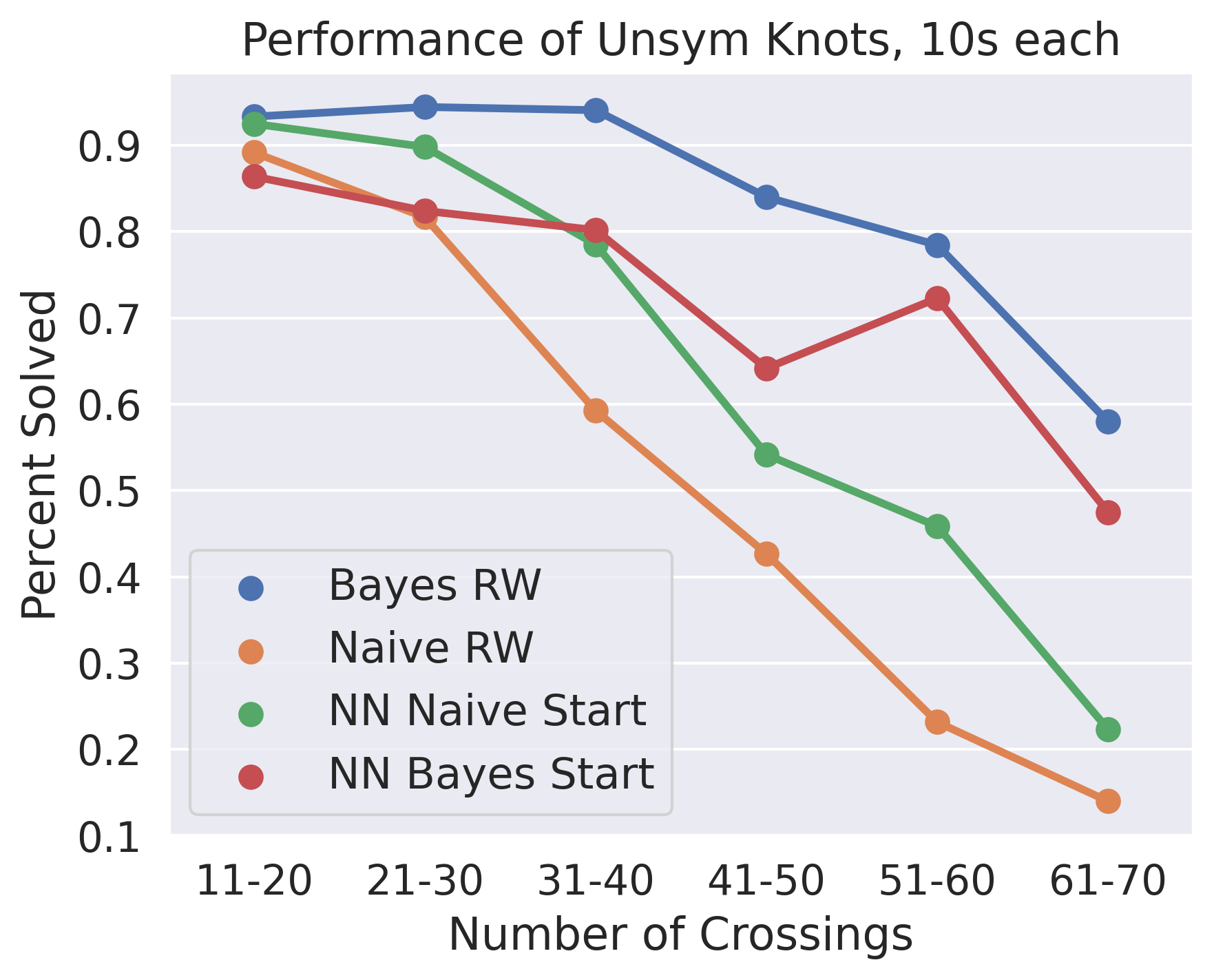}
   \includegraphics[width=.46\textwidth]{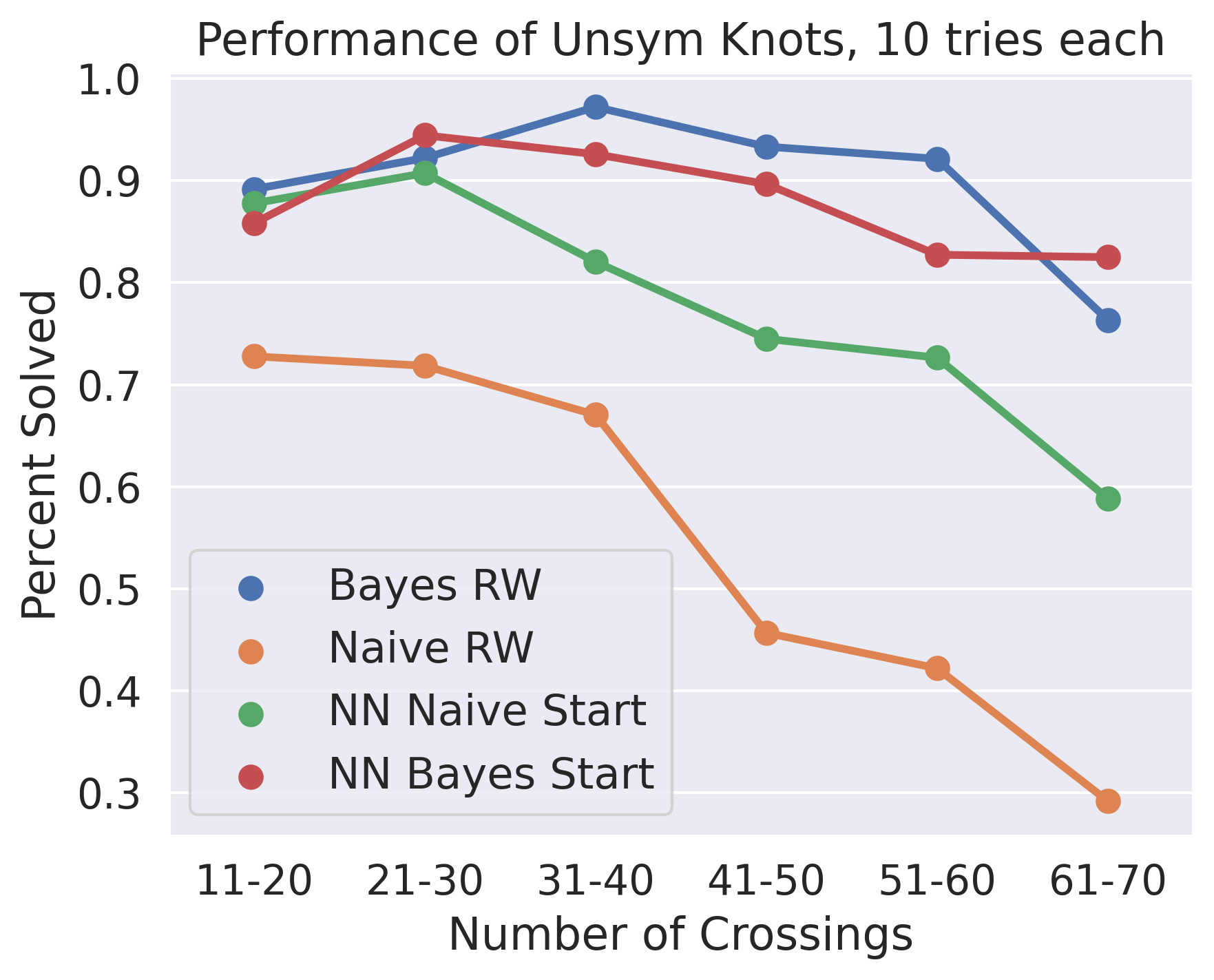}
   \caption{Dependence of performance on the number of crossings for Unsym knots.}
   \label{fig:n_dep_unsym}
\end{figure}

\begin{figure}[!t]
   \centering
   \includegraphics[width=.45\textwidth]{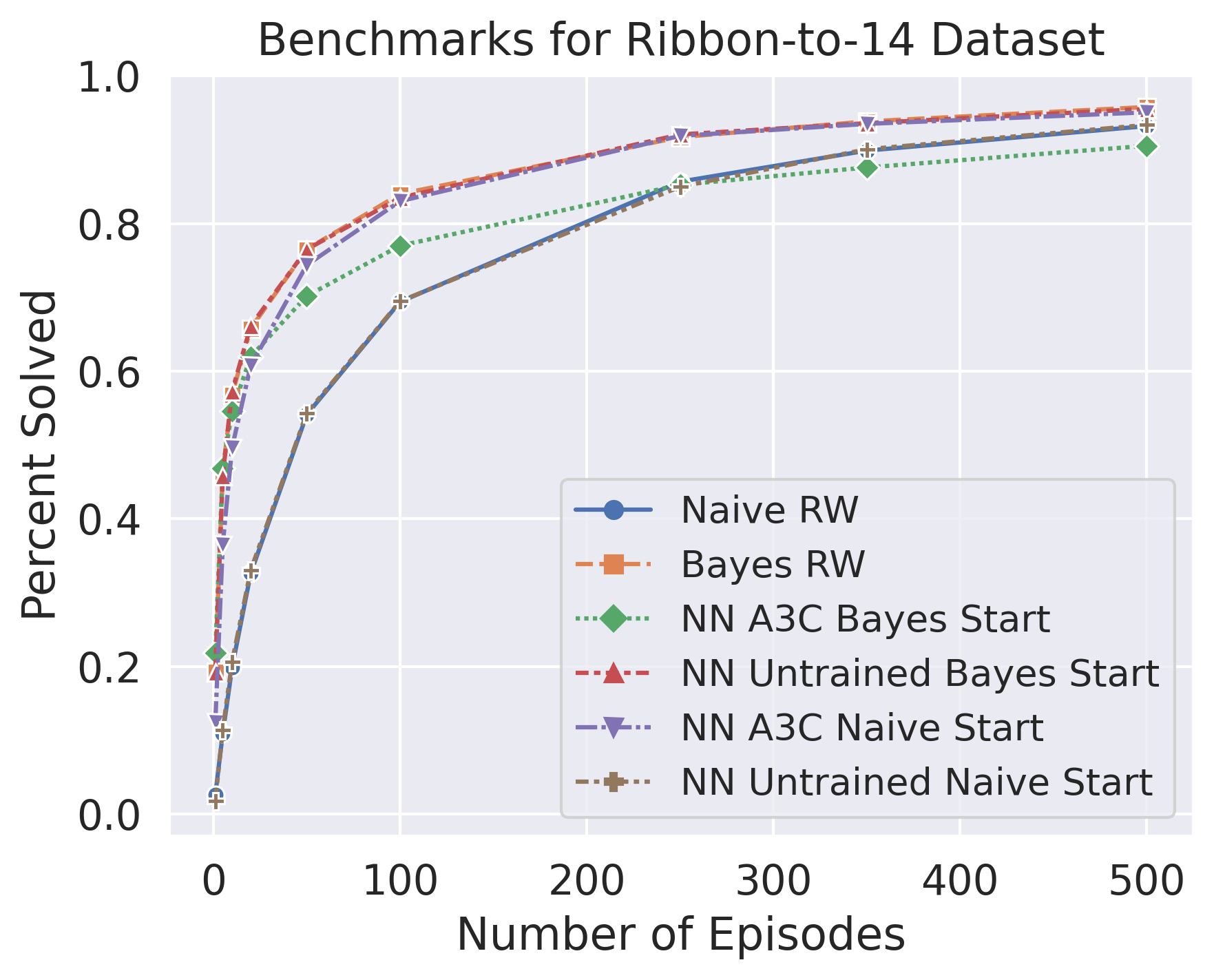}~
   \includegraphics[width=.45\textwidth]{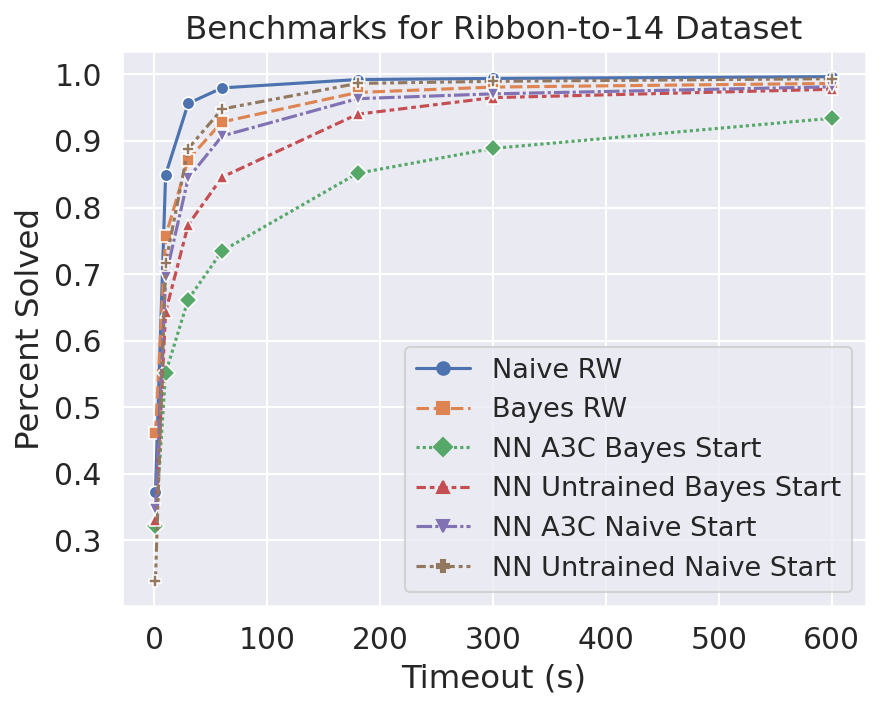}\\
   \includegraphics[width=.45\textwidth]{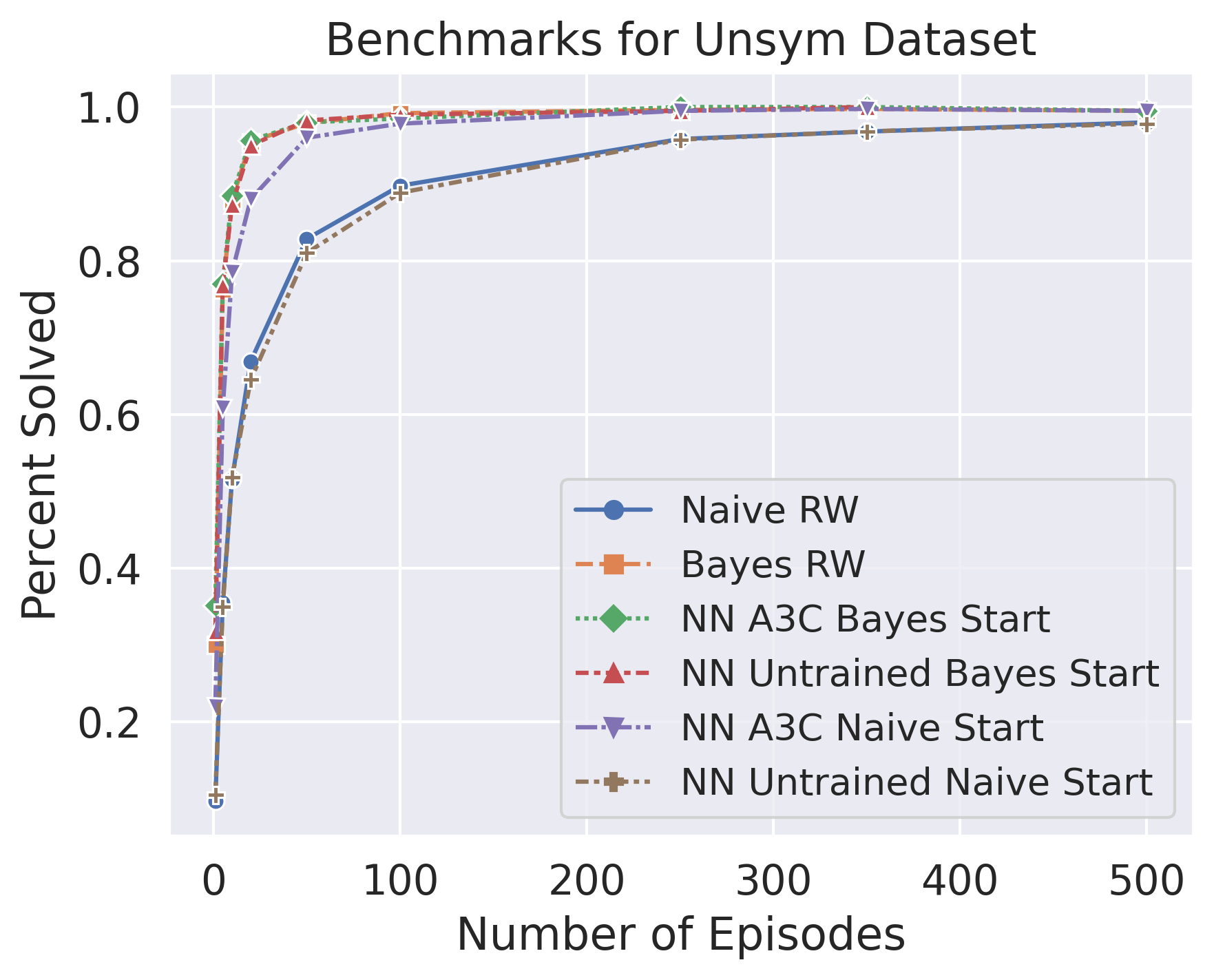}~
   \includegraphics[width=.45\textwidth]{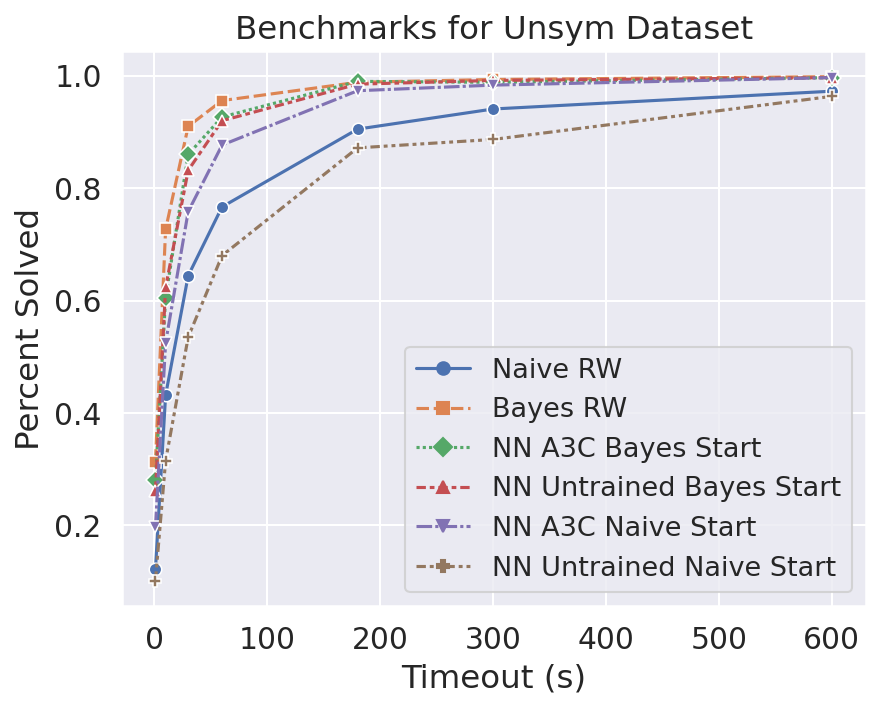}\\
   \includegraphics[width=.45\textwidth]{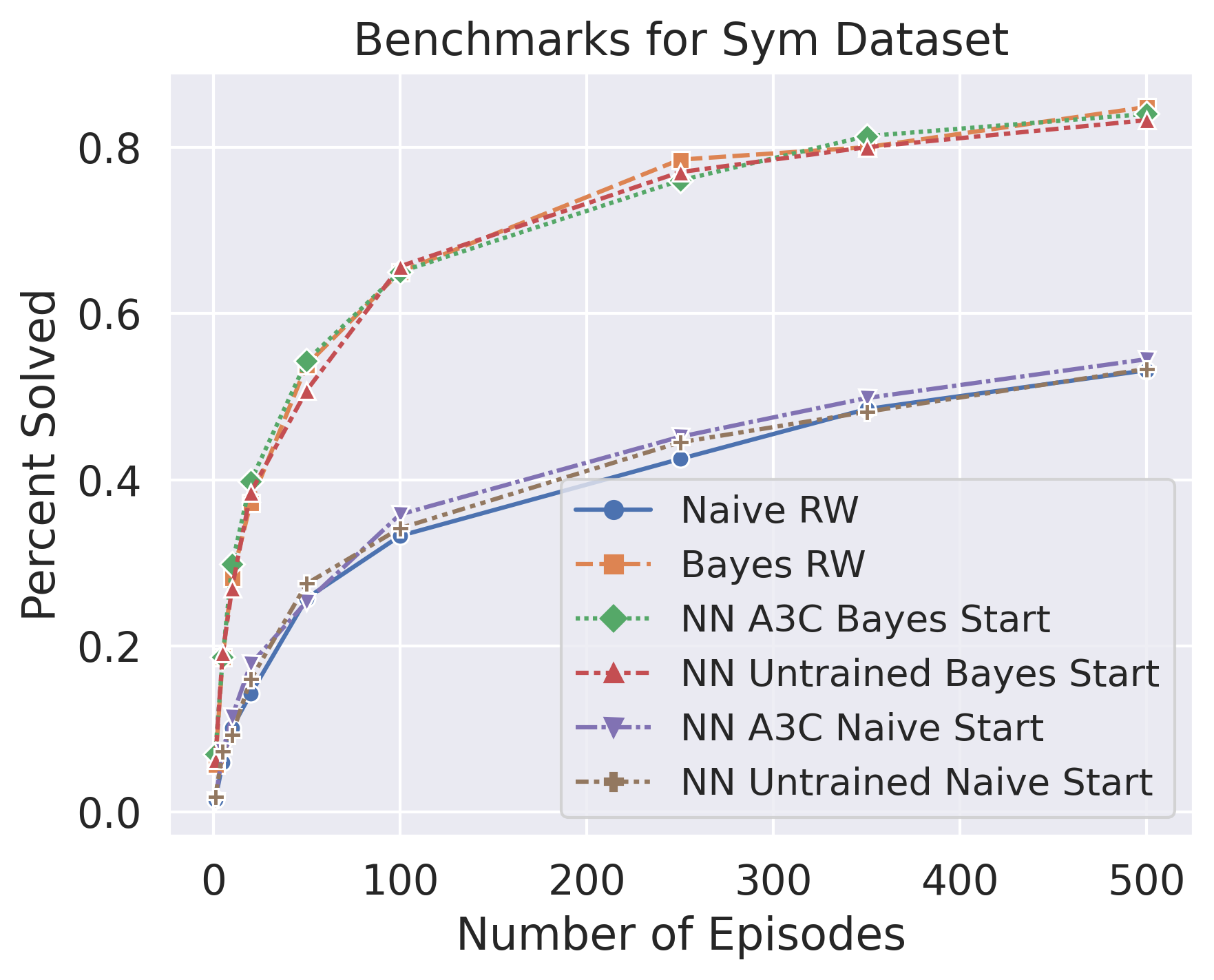}~
   \includegraphics[width=.45\textwidth]{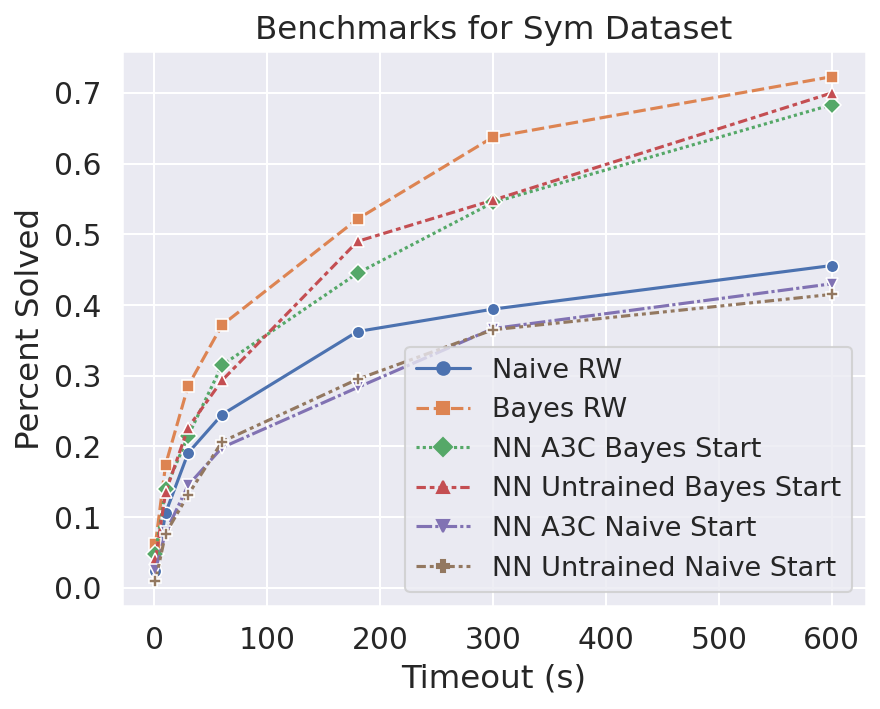}
   \caption{Comparison of six different agents for each knot dataset and benchmark type. The orange line almost coincides with the red line in the plots in the left column.}
   \label{fig:detailed}
\end{figure}

We performed systematic benchmark runs to test the performance of various agents across the different knot datasets. Results are presented in Figures \ref{fig:sym_vs_rto14} - \ref{fig:detailed}.

In much of the following analysis, we will be comparing knot datasets containing knots of different crossing numbers. One of them, Ribbon-to-14, contains all known prime ribbon knots up to 14 crossings, while the other datasets, Sym and Unsym, contain a sample of ribbon knots up to 70 crossings, generated using the techniques outlined in Section~\ref{sec:Generation}. Since the number of actions grows with the number of crossings, we expect higher-crossing knots to be more difficult to handle than lower crossing knots. Hence, we want to start our analysis with a dataset comparison of Ribbon-to-14 and Sym, where we generated around 500 inequivalent (but not necessarily prime) ribbon knots with 6 to 14 crossings; this is to be contrasted with the Sym dataset used throughout, which has 15-70 crossing knots. 

We test the Naive and the Bayes RW on both datasets and time them to measure their performance. We find that both RWs can find ribbon bands for either data set very fast; often in less than one second on a standard MacBook. The times are summarized in the boxplot of Figure~\ref{fig:sym_vs_rto14}, where the box indicates the quartiles of the time distribution while the whiskers show the rest of the distribution and \textsf{x} marks outliers. From the plot, we also see that there does not seem to be a performance difference between all slice knots up to 14 crossings and the ones from the Sym generator, indicating that either representation of the knots is ``equally hard''. Towards the higher crossing end of the datasets, we see the Bayesian optimized RW perform better than the Naive RW. This is likely because it was optimized on knots between 15-70 crossings.

For the rest of this section, when referring to the Sym and Unsym knots, we will mean the datasets containing knots between 15-70 crossings.

In Figure \ref{fig:envelopes} we present the performance of Bayes RW and Naive RW across the different types of knots, with shaded regions to aid visualization and emphasize performance gaps:
\begin{itemize}
\item
In the number-of-episodes benchmarks we clearly see that Unsym knots are easiest to solve, followed by Ribbon-to-14 knots and Sym knots. (This is in spite of the fact that Unsym knots have a much larger number of crossings than Ribbon-to-14 knots.) The performance difference between agents is largest for Sym knots. In all cases the Bayes RW outperforms the Naive RW, despite the fact that the Bayes RW was optimized only for Unsym and Sym knots using a timeout cutoff; \smallskip

\item
In the timeout benchmarks we see that Ribbon-to-14 knots are slightly easier to solve than Unsym knots, and that Sym knots are much more difficult to solve; \smallskip

\item Also in the timeout benchmarks, we see that the Naive RW actually outperforms the Bayes RW on Ribbon-to-14 knots, which can be explained by the fact that the Bayes RW was only optimized on the Unsym and Sym knots.
\end{itemize}

In Figures \ref{fig:n_dep_sym} and \ref{fig:n_dep_unsym} we present the performance of various types of agents on the number of crossings, for Sym knots and Unsym knots, respectively. The agents include Bayes RW, Naive RW, as well as RL agents trained from Bayes and Naive starting distributions, and they are run for a fixed number-of-episodes cutoff or timeout cutoff. The benchmarks are run on the $200$ test-set knots from the respective datasets. We find that:
\begin{itemize}
\item
The Sym knot benchmarks all have decreasing performance with increasing number of crossings, with the Bayes RW and Bayes start RL agent performing best. A concrete takeaway is that to get above $50\%$ accuracy at large crossing number, Sym knots may require running even a Bayes RW for over five minutes. The performance of Bayes (Naive) RW and the RL agent with Bayes (Naive) start are roughly comparable to each other within statistical uncertainty in these small $200$ knot test-set benchmarks; \smallskip
\item
For Unsym knots, however, we see stronger performance correlations: the Naive start RL agent systematically beats the Naive RW, but the Bayes start RL agent is not as strong as the Bayes RW. This suggests that whatever the RL agent is learning when it is given a Naive start, it is not as important as the state-independent distributional shift associated with the Bayes RW optimization, i.e., the relative frequency of the different types of moves is very significant. Notably, the RL agent trained from Naive start significantly outperforms the Naive RW: when given 10 tries per knot the trained agent solved 30\% more Unsym knots with $40$ or more crossings, and it also performs better in the timeout runs, despite the time cost of selecting an action with the neural network.
\end{itemize}

In Figure \ref{fig:detailed} we present six plots, which differ from one another by the choice of one of three knot datasets and whether they are number-of-episodes benchmarks or timeout benchmarks. These are our most detailed plots, as each has six different agent types: a Naive and Bayes RW, a trained and untrained neural network (NN) with Bayes start, and a trained and untrained NN with Naive start. (We sometimes refer to the trained NN as the RL agent.) In the number-of-episodes benchmarks, which make up the left hand side of the figure, we see that:

\begin{itemize}
\item
The curves associated to the Bayes (Naive) RW and the untrained NN with Bayes (Naive) start essentially overlap. This is expected, as the untrained NNs and their associated random walkers have the same probability distribution over actions; the reason for plotting both will become clear when we discuss the timeout benchmarks; \smallskip

\item In the Sym knot dataset,  the trained RL agents with Bayes and Naive start are near the associated RW and Bayes curves, giving evidence that the RL agents have not learned much for Sym knots;
\smallskip
\item
By contrast, for Unsym knots  the RL agent with Naive start significantly outperforms the associated RW and untrained RL agent;
\smallskip
\item
 For Unsym knots, the RL agent with Bayesian start performs similarly to the associated RW and untrained RL agent, suggesting that not much has been learned beyond that of the Bayes RW; \smallskip
\item
Practically all of the Unsym knots are solved by Bayes RW and Naive RW in less than 100s; 
\smallskip
\item
 In the Ribbon-to-14 dataset we see that the RL agent trained from Naive Start performs well, comparable to the Bayes RW and associated untrained RL agent.
 \end{itemize}
 
 In the timeout benchmarks, which make up the right hand side of Figure \ref{fig:detailed}, we see that:
 \begin{itemize}
 \item
There is markedly different performance between the untrained NNs and their associated random walkers, despite having the same probability distributions over actions, with the untrained NN always performing worse. The reason for this is that action selection with the neural network takes longer than that of the Bayes or Naive RW, which does not utilize a neural network. This difference can only show up on timeout benchmarks;
\smallskip
\item
The Bayes RW performs best on the Sym and Unsym datasets, whereas the Naive RW performs best for Ribbon-to-14, which is possible since the Bayes RW is optimized only on the Sym and Unsym datasets;

\smallskip
\item Despite the fact that action selection with the NN introduces a cost that decreases performance, in some cases we see that RL still increases performance over an untrained NN or RW; e.g., for Unsym knots the trained NN with Naive start significantly outperforms the Naive RW and untrained NN with Naive start.
\end{itemize}

In addition to the full list of knots up to 14 crossings and our synthetic datasets ``Sym'' and ``Unsym'', we also tested the algorithm on all alternating knots up to including 20 crossings that are found to be slice using the techniques of ~\cite{Owens:2021aaa}.\footnote{The list is available at \url{https://cat.middlebury.edu/~mathanimations/klo/ribbondisks/}.}  The Bayesian optimized random walker with a 5 min timeout limit identifies about 90\% of these knots as ribbon, with average performance ranging from almost 100\% for alternating knots with up to 14 crossings to around 84\% for alternating knots with 20 crossings. 

The authors of ~\cite{Owens:2021aaa} also list 7 bounty knots for which no slice obstruction could be identified, but for which no ribbon disk was found with their techniques. Our code did not find bands for these 7 knots either.

\section{An RBG family}
\label{sec:RBG}
Recall from Section~\ref{sec:ribbon} that one could disprove SPC4 by finding pairs of knots with the same 0-surgery, such that one is slice and the other is not. In~\cite{MP21}, this strategy was pursued on a family of 3375 pairs of knots (coming from RBG links), with the knots in each pair having the same 0-surgery. The idea suggested in \cite{MP21}, which was inspired from previous work of Freedman, Gompf, Morrison and Walker \cite{MR2657647}, was to use Rasmussen's $s$-invariant from \cite{Rasmussen} to obstruct the sliceness of knots in this family whose companions are slice. However, Nakamura \cite{Nakamura} later showed that the $s$-invariant is not helpful for this purpose. Nevertheless, the possibility remains open for other invariants, so it is worth determining which knots from this family are slice. In principle, a pair of knots with the same $0$-surgery for which the slice status is unknown can be viewed as a potential counterexample to SPC4.

From the RBG family in \cite{MP21}, 2522  pairs can be shown to consist of non-slice (and hence non-ribbon)  knots using various algebraic obstructions. From the remaining 853 pairs, our Bayesian-optimized random walker found 843 of them to consist of ribbon knots. In 5 other pairs, the programs found one knot to be ribbon; the other knot in the pair was then shown to be ribbon using different methods. This left only 5 pairs, and the status of those 10 knots remains unknown. 

\begin{figure}
\centering
\subcaptionbox{$K_G(0,1,-1, -1,1, 0)$}{
\includegraphics[width=.2\textwidth]{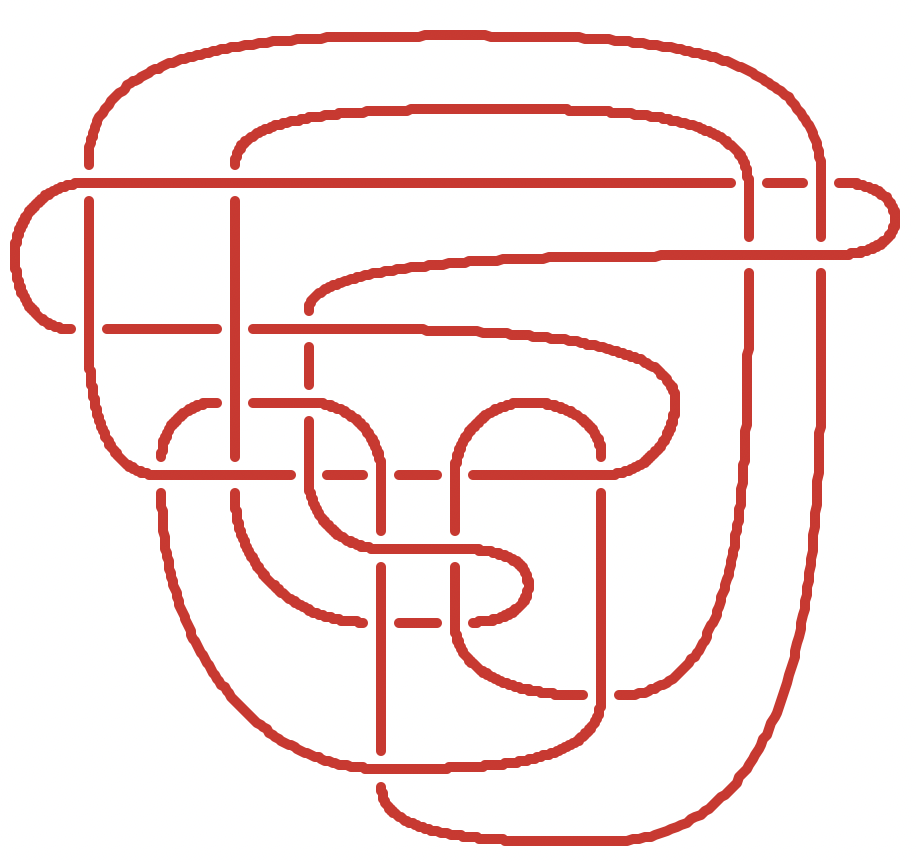}
}\hspace{4mm}
\subcaptionbox{$K_B(0, 1, 2, 0, -1,-1)$}{
\includegraphics[width=.2\textwidth]{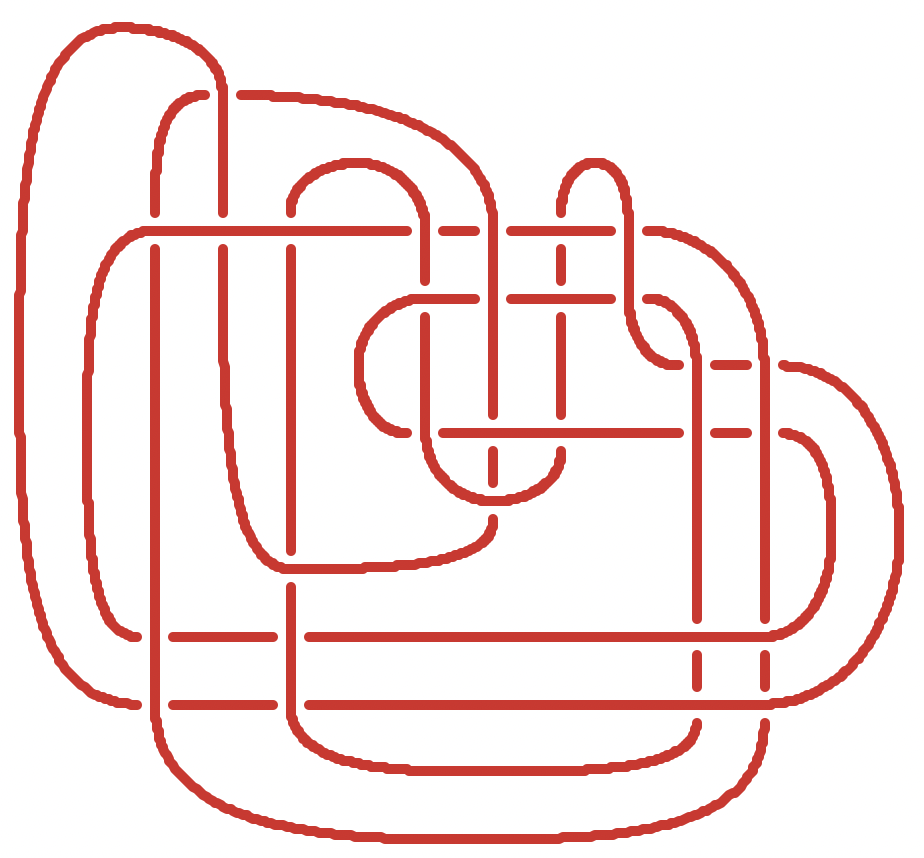}
}\hspace{4mm}
\subcaptionbox{$K_B(0, 0, 2, 0, 0, -1)$}{
\includegraphics[width=.2\textwidth]{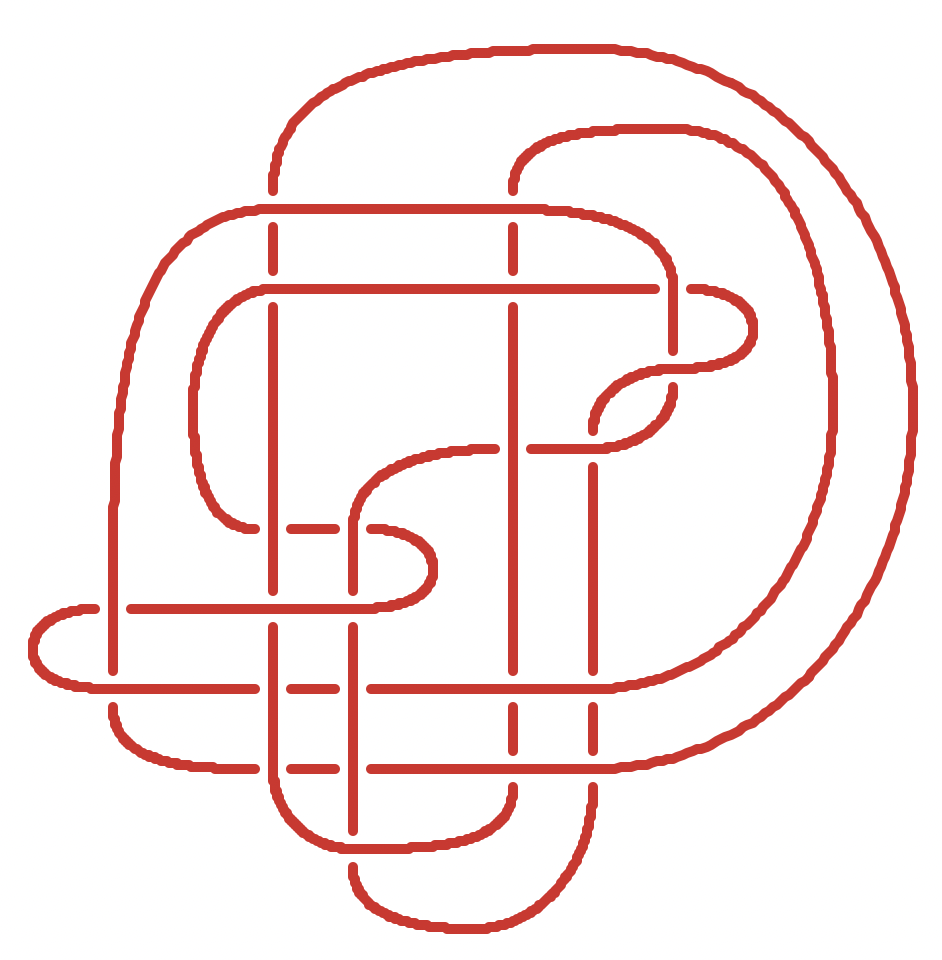}
}\hspace{4mm}
\subcaptionbox{$K_G(2,0,0,-1,2, -1)$}{
\includegraphics[width=.2\textwidth]{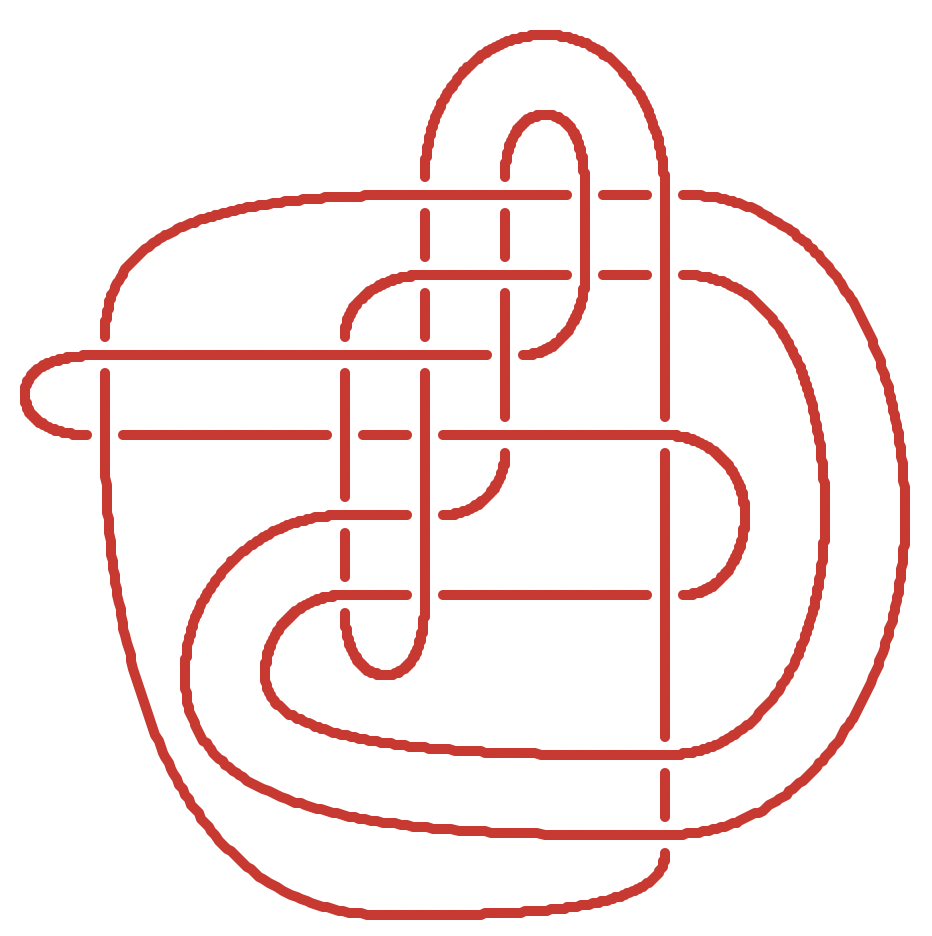}
}
\caption{Ribbon knots not detected by our programs.}
\label{fig:Challenges}
\end{figure}

For completeness, let us discuss how the five knots for which our programs could not find ribbon bands were eventually shown to be ribbon. One of the knots in the list appeared twice, so there were actually only four knots that needed to be analyzed. These four knots are pictured in Figure~\ref{fig:Challenges}, in the notation from \cite{MP21}. They have $23$, $29$, $23$ and $25$ crossings, respectively. 

\captionsetup[sub]{font=small, labelformat=empty}

The first knot, $K_G(0, 1, -1, -1, 1, 0)$, was found to be ribbon (with 3 bands) using an exhaustive search through minimal paths between different segments of the diagram, with a modification of the computer program developed  by Dunfield and Gong \cite{DG}.

The other three knots were shown to be ribbon using an argument suggested to us by Lisa Piccirillo, based on Proposition~\ref{prop:GST}. If $K$ is one of the three given knots, it shares a $0$-surgery with another knot $K'$, which our program found to be ribbon (in all cases, with 2 bands). Hence, one can construct a homotopy four--sphere $W$ by gluing the complement of the slice disk for $K'$ to the trace of the $0$-surgery on $K$, as explained in Section~\ref{sec:ribbon}. The four-dimensional manifold $W$ admits a handle decomposition with three $1$-handles and three $2$-handles. By turning this decomposition upside down, we get one with three $2$-handles and three $3$-handles. One can cancel a $2$-handle against a $3$-handle; afterwards, the attaching link for the remaining $2$-handles is an $R$-link. One of its components is the unknot and the other is the knot $K$, so we can apply Proposition~\ref{prop:GST} to deduce that $K$ is ribbon.

We remark that Proposition~\ref{prop:GST} does not give an easy method of finding the ribbon bands. It remains an interesting challenge to develop a computer program that can find ribbon bands for the three knots $K_B(0, 1, 2, 0, -1,-1)$, $K_B(0, 0, 2, 0, 0, -1)$ and $K_G(2,0,0,-1,2, -1)$. The same goes for the R-links $L_{1,1}$ and $L_{2,1}$ from \cite{Gompf:2010aaa}, which are shown to be ribbon using similar methods to the above.

Thus, every time a knot in a pair was proved to be ribbon, so was its companion; no exotic 4-dimensional spheres were found. While this may seem a negative result, it shows that ML can eliminate many cases, and thus help researchers restrict attention to the remaining few where perhaps a hard-to-find example is lurking. 

In the case at hand, the remaining five pairs whose slice (and ribbon) status is unknown are:
$$ K_{B/G}(0, 0, 0, 1, 2, -1), \ K_{B/G}(0, 0, 0, -1, 2, 1), \ K_{B/G}(0, 0, -2, 0, 0, 1),$$
$$K_{B/G}(-2, 0, 0, -1, 2, -1), \ K_{B/G}(-1, 0, -1, -1, 2, -1).$$ 

The first three of these pairs have $r=0$ in the notation of \cite{MP21}, so by Lemma 5.1(b) in \cite{MP21} the two knots in each pair share the same trace. Therefore, a manifold $W = E(\Delta) \cup (-X(K_1))$ as in Section~\ref{sec:ribbon} would simply be $E(\Delta) \cup (-X(K_2)) = S^4$ rather than an exotic $4$-sphere.

We deduce that, from the original family of 3375 pairs, only two could still potentially produce counterexamples to SPC4: $K_{B/G}(-2, 0, 0, -1, 2, -1)$ and $K_{B/G}(-1, 0, -1, -1, 2, -1)$. (Of course, based on the rest of our analysis, it is unlikely that they do.)

In principle, while the three other pairs 
$$K_{B/G}(0, 0, 0, 1, 2, -1), \ K_{B/G}(0, 0, 0, -1, 2, 1), \ K_{B/G}(0, 0, -2, 0, 0, 1)$$ cannot produce counterexamples to SPC4, they might produce counterexamples to the Slice-Ribbon Conjecture~\ref{conj:slice-ribbon}. Indeed, supposing one of the knots $K_1$ in such a pair is found to be ribbon with a slice disk $\Delta$, the decomposition
$$S^4 = E(\Delta) \cup (-X(K_1)) = E(\Delta) \cup (-X(K_2))$$
would show that $K_2$ bounds an embedded disk in $B^4$, and is therefore slice (but it may not be ribbon). This is similar to the strategy for finding counterexamples to the Slice-Ribbon Conjecture pursued in \cite{MR2740649}. 

We note that many families of knot pairs from RBG links could be studied in the same way, beyond the one considered in \cite{MP21}.

\section{Conclusion}
\label{sec:Conclusions}

We analyzed the performance of six different types of agents across three different datasets of ribbon knots: Ribbon-to-14, Sym, and Unsym. We utilized two different types of benchmarks, in which the agent tries to solve a fixed knot up to some fixed number of episodes or seconds. 

If there is a single takeaway from our analysis, it is that the Bayesian RW performed excellently on all the datasets:
\begin{itemize}
\item  It solved $>99\%$ of the Unsym knots in under $100$ episodes or $5$ minutes;
\item Sym knots were systematically harder, but the Bayes RW still solved $\sim 80\%$ of Sym knots in under $500$ episodes, and $\sim 70\%$ in under $10$ minutes;
\item It achieved a 100\% success rate when given a maximum time cutoff of 10 minutes per knot on the Ribbon-to-14 data set; in fact, the vast majority of these knots can be shown to be ribbon in less than 1~s, as can be seen in Figure~\ref{fig:sym_vs_rto14}. Thus, we recovered the results of \cite{DG} that those 1705 knots with up to 14 crossings are ribbon. 
\end{itemize}

A second takeaway from our analysis is that the data distribution that generates the ribbon knots matters: Sym is much harder to solve than Unsym. This can be seen from the performance of the Bayesian RW noted above, and it can also be seen for other agents. For example, the RL agent trained from a Naive start performed nearly as well as the Bayes RW on Unsym knots in number-of-episodes benchmarks. (The neural network sampling cost slightly increased the performance gap between these two agents in timeout benchmarks). On the other hand, RL did not lead to any improved performance on Sym knots, from neither a Bayes nor a Naive start. The data-dependent failure or success of some techniques over others emphasizes the importance of trying multiple techniques to optimize the likelihood of obtaining new mathematical results.

While our methods were successful for detecting many ribbon knots, interesting challenges remain for future work. In particular, there are knots that we know to be ribbon, but are not recognized as such by our programs. These include, for example:
\begin{itemize}
\item the knots in Figure~\ref{fig:Challenges} that come from the RBG family in \cite{MP21};
\item some of the GST examples in~\cite{Gompf:2010aaa}. There, they exhibit a family of slice links denoted $L_{n, 1}$ that can serve as potential counterexamples to the slice-ribbon conjecture. The first two links in the GST family ($L_{1,1}$ with 18 crossings and $L_{2,1}$ with 40 crossings) were already known to be ribbon (see ~\cite{Gompf:2010aaa}), but unfortunately our algorithms could not find bands to prove they are ribbon. 
\end{itemize}

There also remain knots and links whose ribbon status is unknown, including:
\begin{itemize}
\item the 10 knots left over from the RBG family discussed in Section~\ref{sec:RBG};
\item 21 prime knots with up to 14 crossings;
\item other GST examples from ~\cite{Gompf:2010aaa}, such as $L_{3,1}$ and the slice knot associated to it that is shown in Figure 2 of ~\cite{Gompf:2010aaa};
\item the ``bounty'' alternating knots  mentioned by Owens and Swenton in~\cite{Owens:2021aaa} ;
\item other famous examples such as the positive Whitehead double of the left-handed trefoil. 
 \end{itemize}
 Our programs were unsuccessful at showing these knots and links are ribbon. Of course, it may well be that they are not ribbon, so there is a complementary challenge of finding new powerful obstructions.

To improve the ribbon detection programs, a natural thing to try is to create a different data structure, allowing for more general bands; see Remark~\ref{rem:more}. One could also explore different network architectures. For operation on the dual graph of a knot, graph neural networks might be better suited. Moreover, they would take into account the symmetries of the knot; in our encoding, equivalent knots with isomorphic dual graphs look different to the neural network.

Finally, an interesting challenge is to find better ways of producing ``random'' ribbon knots with large  crossing number. Our Sym and Unsym methods are probably biased towards certain kinds of ribbon knots. Note, for example, that the performance of the programs is much better on the Unsym than on the Sym data sets, indicating that the Unsym knots can be solved with simpler bands. Potentially, machine learning could be used to create new generative models for ribbon knots. For instance, one could use slice obstructions together with our ribbon verifier to perform simulation-based inference. It would use the obstructions or ribbon certificates for each element in an ensemble of knots drawn from a fixed prior to compute a Bayesian posterior conditioned on being ribbon. Samples from the posterior would be more likely to be ribbon, and in general be unrelated to the ribbon distributions (Sym, Unsym, and Ribbon-to-14) utilized in this work.

\appendix
\section{Algorithms}
\label{app:Algorithms}
\begin{algorithm}[t]
\caption{Band Addition MDP Episode}
\label{alg:MDP}
\begin{algorithmic}[1]
\State Initialize state $(G, C, B, T)$ from initial knot diagram
\State $steps \gets 0$
\State Set reward $R \gets 0$
\While{not terminal and $steps < \text{max\_steps}$}
    \State Mask illegal actions based on current state $s_t$
    \State Sample action $a_t \in A$ according to current policy $\pi(a_t \mid s_t)$
    \If{$a_t = \texttt{start}$}
        \State Set starting arc for new band
    \ElsIf{$a_t = \texttt{over}$ \textbf{or} $a_t = \texttt{under}$}
        \State Choose next arc and route band over or under
        \State Update $B$ to include new arc
    \ElsIf{$a_t = \texttt{twist}$}
        \State Insert positive or negative twist; update $T$
    \ElsIf{$a_t = \texttt{end}$}
        \If{band satisfies legality constraints}
            \State Attach band to arc and complete band addition
            \State Simplify link using SnapPy, which applies Reidemeister moves
            \State Update $G$, $C$ accordingly, reset $B$, $T$
            \State $R \gets -(\text{number of crossings of simplified link})$
        \Else
            \State \textbf{terminate episode:} illegal end move, reset $B$, $T$
        \EndIf
    \EndIf
    \If{resulting link is unlink}
        \State \textbf{terminate episode:} win
    \EndIf
    \If{hyperparameter constraints violated}
        \State \textbf{terminate episode:} loss
    \EndIf
    \State $steps \gets steps + 1$
\EndWhile
\State Return reward $R$
\end{algorithmic}
\end{algorithm}

\clearpage

\bibliography{refs}
\bibliographystyle{amsalpha}
\addresseshere

\end{document}